\input epsfx.tex
\overfullrule=0pt


\magnification=1400
\hsize=12.5cm
\vsize=16.5cm
\hoffset=-0.3cm   
\voffset=1.0cm    

\baselineskip 16 true pt

\font \smcaps=cmbx10 at 12 pt


\font \caps=cmbx10 scaled 1200

\font \smcaps=cmcsc10 
\font \pecaps=cmcsc10 at 9 pt

\newtoks \hautpagegauche  \hautpagegauche={\hfil}
\newtoks \hautpagedroite  \hautpagedroite={\hfil}
\newtoks \titregauche     \titregauche={\hfil}
\newtoks \titredroite     \titredroite={\hfil}
\newif \iftoppage         \toppagefalse   
\newif \ifbotpage         \botpagefalse    
\titregauche={\pecaps   R\'egis Msallam and Fran\c{c}ois Dubois  } 
\titredroite={\pecaps   Model for coupling  a   perfect flow   with an
acoustic boundary layer   } 
\hautpagegauche = { \hfill \the \titregauche  \hfill  }
\hautpagedroite = { \hfill \the \titredroite  \hfill  }
\headline={ \vbox  { \line {  
\iftoppage    \ifodd  \pageno \the \hautpagedroite  \else \the
\hautpagegauche \fi \fi }     \bigskip  \bigskip  }}
\footline={ \vbox  {   \bigskip  \bigskip \line {  \ifbotpage  
\hfil {\oldstyle \folio} \hfil  \fi }}}
 



\def\R{{\rm I}\! {\rm R}}


\def\sqr#1#2{{\vcenter{\vbox{\hrule height.#2pt \hbox{\vrule 
width .#2pt height#1pt \kern#1pt \vrule width.#2pt}
\hrule height.#2pt}}}}

\def\mod#1{\setbox1=\hbox{\kern 3pt{#1}\kern 3pt}%
\dimen1=\ht1 \advance\dimen1 by 0.1pt \dimen2=\dp1 \advance\dimen2 by 0.1pt
\setbox1=\hbox{\vrule height\dimen1 depth\dimen2\box1\vrule}%
\advance\dimen1 by .1pt \ht1=\dimen1
\advance \dimen2 by .01pt \dp1=\dimen2 \box1 \relax}

\def\nor#1{\setbox1=\hbox{\kern 3pt{#1}\kern 3pt}%
\dimen1=\ht1 \advance\dimen1 by 0.1pt \dimen2=\dp1 \advance\dimen2 by 0.1pt
\setbox1=\hbox{\kern 1pt  \vrule \kern 2pt \vrule height\dimen1 depth\dimen2\box1
\vrule
\kern 2pt \vrule \kern 1pt  }%
\advance\dimen1 by .1pt \ht1=\dimen1
\advance \dimen2 by .01pt \dp1=\dimen2 \box1 \relax}

\font \caps=cmbx10 at 14 pt
\font \smcaps=cmbx10 at 12 pt

\def\br {\break}

  \def \page #1{\unskip\leaders\hbox to 1.3 mm {\hss.\hss}\hfill $\,$ {\oldstyle #1}}

\rm 


$~$

\bigskip
\bigskip
\bigskip
\bigskip
\centerline{\caps  Mathematical model for coupling } \medskip
\centerline{\caps  a quasi-unidimensional perfect flow } \medskip
\centerline{\caps  with an acoustic boundary layer}

\bigskip
\bigskip
\bigskip

\centerline  {\smcaps    R\'egis Msallam \footnote{$ ^{^{\displaystyle \diamond}}$}
{\rm Institut de Recherche et de Coordination  Acoustique /  Musique,
\smallskip  \vskip -2pt 1 place Stravinsky, F-75004 $\,$ Paris, France ; 
Regis.Msallam@ircam.fr}  $\,$ and $\,$  Fran\c{c}ois Dubois 
\footnote{$ ^{^{\displaystyle \smcaps \ast}}$}
{\rm Applications Scientifiques du Calcul Intensif,   b\^at. 506, BP 167,  F-91~403 
\smallskip  \vskip -3pt 
Orsay  Cedex, Union Europ\'eenne ;   $\,\,$   dubois@asci.fr. }}

\bigskip
\centerline{July 26, 1999, revised version  June 12, 2002. 
\footnote{$ ^{^{\displaystyle \dag}}$}
{\rm  Rapport  n$^{\rm o}$ 326-99 de l'Institut A\'ero Technique du 
Conservatoire National 
\smallskip  \vskip -3pt 
des Arts et M\'etiers \`a  Saint Cyr l'Ecole.
  Edition du 10   juin 2010.}}

\null\vskip 1cm

\bigskip \bigskip \noindent {\smcaps Abstract}

Nonlinear acoustics of wind instruments conducts to study unidimensional
fluid flows. From physically relevant approximations that are modelized with  the thin layer
Navier Stokes equations, we propose a coupled model where perfect fluid flow is described
by the Euler equations of gas dynamics and viscous and thermal boundary layer
is modelized by a  linear equation. We describe numerical discretization, validate the
associated software by comparison with analytical solutions and consider musical
application of strongly nonlinear waves in the trombone. 

\bigskip 
\vfill \eject 
\noindent {\smcaps R\'esum\'e}

L'acoustique non lin\'eaire des instruments \`a vents conduit \`a  \'etudier
les \'ecou\-lements filaires monodimensionnels. A partir d'approximations physiquement\br
r\'ealistes qui sont prises en compte par les  \'equations de Navier Stokes de couche mince, nous 
proposons   un mod\`ele coupl\'e o\`u le fluide parfait est d\'ecrit par les \'equations
d'Euler de la dynamique des gaz et le fluide visqueux  et conducteur de chaleur par une
\'equation lin\'eaire de couche limite. Nous d\'etaillons la discr\'etisation num\'erique
retenue et validons le logiciel d\'evelopp\'e gr\^ace \`a des solutions analytiques avant
d'aborder l'application musicale
\`a la propagation d'ondes fortement non lin\'eaires dans le trombone.

\toppagetrue  
\botpagetrue    

\bigskip
\bigskip \bigskip \noindent {\bf Key words :} fluid mechanics, nonlinear
acoustics, Euler equations, boundary layer, finite differences.


\bigskip 
\bigskip \noindent {\smcaps Contents}   \smallskip  

1) Introduction \page {$\,\,$2}

2) Thin Layer Navier Stokes equations \page {$\,\,$3}

3)  Perfect fluid for main flow  \page {$\,\,$10}

4)  Acoustic boundary layer   \page {$\,\,$13}

5)  The coupled problem  \page {$\,\,$14}

6)   Generalization to axisymmetric geometry   \page {$\,\,$17}

7)   Numerical approximation of the coupled problem   \page {$\,\,$20}

8)   First test cases   \page {$\,\,$27}

9)   Conclusion, aknowledgments   \page {$\,\,$35}

10)  References   \page {$\,\,$35}


\bigskip \bigskip \bigskip
\noindent {\smcaps 1) $\quad$ Introduction.} 
\smallskip \noindent $\bullet \quad$
In this paper, we study simple models of nonlinear acoustic flows in  cylindric or
axisymmetric  ducts. Our objective is to take into account several physical effects such
compressibility of the air, viscous  dissipation and thermal conduction, expecially
in the vicinity of the wall. We consider the flow of a newtonian compressible fluid
in a two-dimensional pipe. In a first approximation, the variation of physical
fields in the transverse direction can be neglected and an appropriate physical
model for such a flow is given by unidimensional equations of gas dynamics (see,
{\it e.g.} Landau-Lifschitz [LL53]). This model is appropriate for the description of
nonlinear waves in shock dynamics (see {\it e.g.} Courant-Friedrichs [CF48]) and also for
weaker waves in nonlinear acoustics (Whitham [Wh74]). Nevertheless, such a model
neglects all phenomena that can appear in the boundary layer.

\bigskip \noindent $\bullet \quad$ 
The boundary layer is the region located near the wall
(at a distance of the order of the boundary layer thickness $\delta$) where
viscous dissipation and thermal conduction have to be considered. The important
role of the boundary layer in duct acoustics has been  studied by Chester [Ch64].
Usually, linearized boundary layer equations are considered and the associated
mathematical model is the heat equation whose solution can
be explicited using a convolution kernel in time. 

\bigskip \noindent $\bullet \quad$
In the domain of acoustics, nonlinear and dissipative
effects are usually taken into account via generalized Burgers equations as
suggested by Blackstock [Bl85] ; this scalar model contains  a source term which
is, in the case of ducts, a convolution kernel giving an explicit solution  of
the linear model of acoustic boundary layer. We refer to Makarov-Ochmann [MO97]
for a review of the fundamental results. 

\bigskip \noindent $\bullet \quad$
We here focus on the fact that the modelling of an
acoustic flow in a pipe can be conducted as a coupling between a perfect fluid
and a boundary layer. We refer to Le Balleur [LB80], Zeytounian [Ze92] and
Aupoix-Brazier-Cousteix [ABC92] for classical approaches developed in the context
of aerodynamics applications. In this paper, we have been  inspired  by
these coupling techniques for pipe flow problem in nonlinear acoustics.

\bigskip
\bigskip
\noindent {\smcaps 2) $\quad$  Thin Layer Navier Stokes equations.} 
\smallskip \noindent $\bullet \quad$ 
We consider the   geometry of a two-dimensional pipe of characteristic longitudinal length
equal to $L$. The thickness of the duct is $2 h$ and our basic hypothesis is that the ratio
$ \displaystyle \, {{h}\over{L}} \, $  is small~: 

\smallskip \noindent  (2.1) $\qquad \displaystyle
{{h}\over{L}} \qquad << \quad 1 \,. $

\smallskip \noindent
We consider also some longitudinal length $\, \Lambda \,$ and some length  $ \, l \, $ in
the  transverse direction ($\, 0  <   l  \leq  h \, $)
that are used for the adimensionalization of the equations (see Figure 1). We distinguish
between the two types of flows associated to  aerodynamics  and acoustics
applications. 

\bigskip 
\centerline { \epsfysize=4cm    \epsfbox  {fig1.epsf} }
\smallskip  \smallskip
\centerline { {\bf Figure 1} \quad Channel with characteristic lenghts. }
\smallskip

\bigskip \noindent $\bullet \quad$
In aerodynamics, we suppose simply~:  

\smallskip \noindent  (2.2) $\qquad \displaystyle
\Lambda \,\, = \,\, L \, . \, $

\smallskip \noindent 
Moreover, the distance  $ \, l \, $ is a distance characteristic of the maximum of the
boundary layer thickness (see {\it e.g.} Schlichting [Sc55])~:  

\smallskip \noindent  (2.3) $\qquad \displaystyle
l = \,5\,\sqrt{{\mu\,L}\over{\rho\,U}}\quad.  $

\smallskip \noindent
In previous expression,  $\rho\,$, $\,\mu\,$ and  $\,U\,$ are respectively  the density, the
viscosity and the amount velocity of the flow and we introduce also the so-called Reynolds
number in aerodynamics~:  

\smallskip \noindent  (2.4) $\qquad \displaystyle
{\it R_{\rm e}^{\rm aero}} \,\, = \,\, {{\rho \, U \, L}\over{\mu}}\, . \,    $

\smallskip \noindent
We note that  for extremely thin pipes, the boundary layer occupies all the duct, that is $
l \simeq h$. In all cases, we suppose that 

\smallskip \noindent  (2.5) $\qquad \displaystyle
\epsilon \,\, \equiv \,\, {{l}\over{\Lambda}} \,\, \approx \,\, {{1}\over{\sqrt{\it R_{\rm
e}^{\rm aero}}}} \qquad << \quad 1 \,. $

\bigskip \noindent $\bullet \quad$
In acoustics, if condition

\smallskip \noindent  (2.6) $\qquad \displaystyle
{{h}\over{ \lambda }} < {{1}\over{4}}  \, \, $

\smallskip \noindent  
is satisfied, the waves propagate only  along the axial direction (Pierce [Pi81], see also
Bruneau [Br98]) and it is natural to consider the length wave $ \, \lambda \,$  as a
reference length for the axial direction. We set~:  

\smallskip \noindent  (2.7) $\qquad \displaystyle
\Lambda \,\, = \,\, \lambda   \, .  \, $

\smallskip \noindent 
On the other hand,  distance $\, l \, $ is the natural length constructed from viscosity
coefficient $\, \mu \, $, density of the air $\, \rho_{0} \, $ and sound celerity
$\, c_{0} \, $  at usual thermodynamic conditions for pressure and temperature~: 
$\,p_{0}=~1\,$atmosphere and $\,T_{0}=  300\,$ Kelvins. We set  

\smallskip \noindent  (2.8) $\qquad \displaystyle
 l = {{\mu}\over{\rho_{0}\,c_{0}}} \,  $

\smallskip \noindent 
and we introduce also the acoustic Reynolds number $ \, {\it R_{\rm e}^{\rm acou}} \, $
defined simirarily to expression (2.4) : 

\smallskip \noindent  (2.9) $\qquad \displaystyle
{\it R_{\rm e}^{\rm acou}} \,\, = \,\,  {{\rho_{0} \,  c_{0}  \, \lambda}\over{\mu}}\, .
\,  $

\smallskip \noindent 
The ratio $\, \displaystyle {{l}\over{\Lambda}} \,$ between right hand sides of
expressions  (2.8) and (2.7) satisfies the following hypothesis : 

\smallskip \noindent  (2.10) $\qquad \displaystyle
\epsilon \,\, \equiv \,\, {{l}\over{\Lambda}} \,\, \approx \,\, {{1}\over{\it R_{\rm
e}^{\rm acou}}} \qquad << \quad 1 \,. $

\smallskip \noindent 
In practice, $\,l\simeq 10^{-8}\,m\,$ and if the frequency of the acoustic
wave (with wave length~$\lambda$) is less than $1 \, {\rm Ghz}$, condition (2.10) is
satisfied. More precisely, we set with Bruneau, Herzog, Kergomard and Polak
[BHKP89], 

\smallskip \noindent  (2.11) $\qquad \displaystyle
  l_{vh} \,\,=\,\, \Bigl( {{4}\over{3}}\mu\, +\mu_{v}\Bigr)
{{1}\over{\rho_{0}\,c_{0}}}\,\, + \,\,(\gamma-1)
{{k}\over{{\rho_{0}\,c_{0}\,C_{p}}}} $

\smallskip \noindent  
where  $\,\mu_{v}\,$ is the volumic viscosity, $\,\gamma\,=\,7/5\,$ is the ratio
of specific heats, $\,k\,$ the thermic conductivity and $\,C_{p}\,$ the
calorific capacity at constant pressure. In the air the
volumic viscosity $\,\mu_{v}\,$ is negligeable compared to viscosity
$\,\mu\,$ and $\,(\gamma-1)\,k\,/\,C_{p}\,$ is of the order of viscosity
$\,\mu\,$ {\it i.e.} the Prandtl number (see {\it e.g.} Schlichting [Sc55]) is of the order
of $\,1\,$. Therefore, 

\smallskip \noindent  (2.12) $\qquad \displaystyle
$$  l \approx l_{vh} . $

\smallskip \noindent 
that enforces hypothesis (2.8). 

\bigskip \noindent $\bullet \quad$
The flow is supposed to satisfy the Navier Stokes equations of conservation of mass,
impulse and energy. Recall that the unknowns are density $\,\rho\,$, velocity
$\,(u,v)$, pressure $\,p\,$ and internal specific energy $\,e\,$. The thermodynamic variables are supposed to satisfy the
state equation for the air that takes the classical form of a perfect gas equation

\smallskip \noindent  (2.13) $\qquad \displaystyle
p\,\,=\,\,(\gamma-1)\,\rho\,e\,.$

\smallskip \noindent
In the following, we neglect the volumic viscosity $\,\mu_{v}\,$ and assume
that the Stokes hypothesis concerning the two viscosities is valid. In
consequence, the classical analytic expression of the Navier-Stokes equations 
({\it e.g.} Landau and Lifschitz [LL53])  takes the form

\smallskip \noindent  (2.14) $\quad \displaystyle
{{\partial\,\rho}\over{\partial\,t}}\,\,+\,\,{{\partial}\over{\partial\,x}}
(\rho\,u)\,\,+\,\,{{\partial}\over{\partial\,y}}(\rho\,v)\,\,=\,\,0 $

\smallskip \noindent  (2.15) $\quad \displaystyle
{{\partial}\over{\partial\,t}}(\rho\,u) +
{{\partial}\over{\partial\,x}}\bigl(\rho\,u^{2}\,+\,p\,\bigr)+
{{\partial}\over{\partial\,y}}\bigl(\rho\,u\,v\bigr) \, = \, \mu\, \biggl[
{{\partial}\over{\partial\,x}}\Bigl({{4}\over{3}}{{\partial\,u}\over{\partial\,x}}
\,+\,{{1}\over{3}}{{\partial\,v}\over{\partial\,y}} \Bigr) + 
{{\partial^{2}\,u}\over{\partial\,y^{2}}}  \,\biggr] $

\smallskip \noindent  (2.16) $\quad \displaystyle
{{\partial}\over{\partial\,t}}(\rho\,v) +
{{\partial}\over{\partial\,x}}\bigl(\rho\,u\,v\,\bigr) + 
{{\partial}\over{\partial\,y}}\bigl(\rho\,v^{2}\,+\,p\,\bigr) \, = \,  \mu\,
\biggl[\,{{\partial^{2}\,v}\over{\partial\,x^{2}}} + 
{{\partial}\over{\partial\,y}}\Bigl({{1}\over{3}}{{\partial\,u}\over{\partial\,x}}
\,+\,{{4}\over{3}}{{\partial\,v}\over{\partial\,y}} \Bigr) \,\biggr] $

\setbox20=\hbox{$\displaystyle 
{{\partial}\over{\partial\,t}} \biggl( \rho\,\Bigl(
e+{{1}\over{2}} \bigl( u^{2}+v^{2}\bigr) \Bigr) \biggr)\,\,+\,\,
{{\partial}\over{\partial\,x}}\biggl(\rho\,u\,\Bigl(
e+{{1}\over{2}} \bigl( u^{2}+v^{2}\bigr) \Bigr)\,+\,p\,u\biggr)\,\,+\qquad  $}
\setbox21=\hbox{$\displaystyle  
 \,\,+ \,\,{{\partial}\over{\partial\,y}}\biggl(\rho\,v\,\Bigl(
e+{{1}\over{2}} \bigl( u^{2}+v^{2}\bigr) \Bigr)\,+\,p\,v\biggr)\,\,=\,\, $}
\setbox22=\hbox{$\displaystyle \qquad \qquad  \qquad \qquad 
\,\,=\,\,  \mu\,{{\partial}\over{\partial\,x}} \,\biggl[\,u\,\biggl(
{{4}\over{3}}{{\partial\,u}\over{\partial\,x}}
\,-\,{{2}\over{3}}{{\partial\,v}\over{\partial\,y}}\biggr)\,+\,
 v\,\biggl( {{\partial\,u}\over{\partial\,y}} \,+\,
{{\partial\,v}\over{\partial\,x}} \biggr) \biggl]   \,+   $}
\setbox23=\hbox{$\displaystyle 
\, + \,  \mu \,{{\partial}\over{\partial\,y}} \, \biggl[ \,u   \biggl(
{{\partial\,v}\over{\partial\,x}}\,+\,{{\partial\,u}\over{\partial\,y}} \biggr)
+ v \biggl({{4}\over{3}} {{\partial\,v}\over{\partial\,y}}
\,-\,{{2}\over{3}}{{\partial\,u}\over{\partial\,x}}\biggr)\, \biggl]  \,  
+ \, k\,\biggl( {{\partial^{2}\,T}\over{\partial\,x^{2}}} \,+\, 
{{\partial^{2}\,T}\over{\partial\,y^{2}}} \biggl)\,\,. $}
\setbox30= \vbox {\halign{#\cr \box20 \cr \box21 \cr \box22 \cr \box23 \cr}}
\setbox31= \hbox{ $\vcenter {\box30} $}
\setbox44=\hbox{\noindent  (2.17) $ \left\{ \!\! \box31 \right. $}  
\smallskip \noindent $ \box44 $

\bigskip \noindent $\bullet \quad$
We detail the way we adimensionalize the Navier Stokes equations. First we have two
length scales $\, \Lambda \,$ and $\, l \,$ for longitudinal and transverse directions
respectively ; we denote by $\,\overline{x}\,$ and $\,\overline{y}\,$ these two space
variables without dimension~: 

\smallskip \noindent
\smallskip \noindent  (2.18) $\qquad \displaystyle
\overline{x}\,\,=\,\,{{x}\over{\Lambda}}$
\smallskip \noindent  (2.19) $\qquad \displaystyle
\overline{y}\,\,=\,\, {{y}\over{l}} \,.$

\smallskip \noindent
Second, we introduce some longitudinal reference velocity $\,U\,$. This velocity
defines a time reference $\,\tau\,$ and an adimensionalized time
$\,\overline{t}\,$ according to

\smallskip \noindent
\smallskip \noindent  (2.20) $\qquad \displaystyle
\tau\,\,=\,\,{{\Lambda}\over{U}} $
\smallskip \noindent  (2.21) $\qquad \displaystyle
\overline{t}\,\,=\,\, {{t}\over{\tau}} \,.$

\smallskip \noindent
If $\,\Lambda =\lambda\,$ is the length wave and $\,U=c_{0}\,$ is a typical choice for
the  adimensionnalization of velocity in acoustics, then $\,\tau\,$ is the period of
the wave,  {\it i.e.} the time for the acoustic perturbation  to travel one length wave.
We introduce a second reference velocity $\,V\,$ associated to this time
$\,\tau\,$  and the transverse distance $\,l\,$:

\smallskip \noindent  (2.22) $\qquad \displaystyle
 V\,\,={{l}\over{\tau}}\,.$

\smallskip \noindent 
If a particle travels distance $\,\Lambda\,$ with axial velocity $\,U\,$ during time
$\,\tau\,$, it travels distance $\,l\,$ with transverse velocity $\,V\,$ during the
same time interval. Therefore, we define dimensionless velocities 
$\,\overline{u}\,$ and $\,\overline{v}\,$ according to

\smallskip \noindent  (2.23) $\qquad \displaystyle
\overline{u}\,\,=\,\,{{u}\over{U}}$
\smallskip \noindent  (2.24) $\qquad \displaystyle
\overline{v}\,\,=\,\,{{v}\over{V}} \,\,=\,\,{{1}\over{\epsilon}}
\,{{v}\over{U}}  $

\smallskip \noindent  with

\smallskip \noindent  (2.25) $\qquad \displaystyle
\epsilon\,\,={{l}\over{\Lambda}}\,. $

\bigskip \noindent $\bullet \quad$  
For the adimensionalization of convective	terms, the reference for density is the
density $\,\rho_{0}\,$ of the air at the usual conditions and reference for
pressure is associated with the dynamic pressure	$\,\rho_{0}\,U^{2}\,$. We set 

\smallskip \noindent  (2.26) $\qquad \displaystyle
\overline{\rho}\,\,=\,\,{{\rho}\over{\rho_{0}}}  $
\smallskip \noindent  (2.27) $\qquad \displaystyle
\overline{p}\,\,=\,\,{{p}\over{\rho_{0}\,U^{2}}}\,. $

\smallskip \noindent
The reference for internal energy is chosen in order to maintain the
validity of the state equation (2.13) after adimensionalization. We set

\smallskip \noindent  (2.28) $\qquad \displaystyle
\overline{e}\,\,={{e}\over{U^{2}}}\, $

\smallskip \noindent 
and we deduce from (2.13), (2.26) (2.27) and (2.28) the state equation between
these new variables~: 

\smallskip \noindent  (2.29) $\qquad \displaystyle
\overline{p}\,\,=\,\, (\gamma-1)\,\overline{\rho}\,\overline{e} \,. $

\bigskip \noindent $\bullet \quad$ 
The Reynolds number $\,{\it R_{\rm e}}\,$ appears from the dissipation terms in the
momentum equations (2.15) and (2.16)

\smallskip \noindent  (2.30) $\qquad \displaystyle
{\it R_{\rm e}}\,\,=\,\,{{\rho_{0}\,U\,\Lambda}\over{\mu}}\,, $

\smallskip \noindent
the Prandtl number $\,{\it P_{\rm r}}\,$ is defined from the heat fluxes in the
energy equation (2.17)

\smallskip \noindent  (2.31) $\qquad \displaystyle
{\it P_{\rm r}}\,\,=\,\,{{\mu\,C_{p}}\over{k}} $

\smallskip \noindent 
and a reference scale for temperature is defined by $\displaystyle
\,{{U^{2}\over{C_{p}}}}\,$~: 

\smallskip \noindent  (2.32) $\qquad \displaystyle
\overline{T}\,\,=\,\,{{C_{p}\,T}\over{U^{2}}}\,.$

\smallskip \noindent
The Joule-Gay Lussac law for polytropic gas can be rewritten in terms of
dimensionless energy  $\,\overline{e}\,$ and temperature $\,\overline{T}\,\,$
according to

\smallskip \noindent  (2.33) $\qquad \displaystyle
\overline{e}\,\,=\,\,\gamma\,  \overline{T}\,. $

\bigskip \noindent $\bullet \quad$
Then the adimensionalized Navier Stokes equations take the following form 

\smallskip \noindent  (2.34) $\qquad \displaystyle
{{\partial\,\overline{\rho}}\over{\partial\,\overline{t}}}\,\,+\,\,
{{\partial}\over{\partial\,\overline{x}}} \,
\bigl( \overline{\rho}\,\overline{u} \bigr) \,\,+\,\,
{{\partial}\over{\partial\,\overline{y}}}\,
\bigl( \overline{\rho}\,\overline{v} \bigr) \,\,=\,\,0  $

\setbox20=\hbox{$\displaystyle 
{{\partial}\over{\partial\,\overline{t}}}(\overline{\rho}\,\overline{u})\,\,+\,\,
{{\partial}\over{\partial\,\overline{x}}}\bigl(\overline{\rho}\,\overline{u}^{2}\,
+\,\overline{p}\,\bigr)\,\,+\,\,
{{\partial}\over{\partial\,\overline{y}}}\bigl(\overline{\rho}\,\overline{u}\,
\overline{v}\bigr)\,\,= $}
\setbox21=\hbox{$\displaystyle   \qquad \qquad  \qquad \qquad \qquad 
\,\,=\,\,{{1}\over{\it R_{\rm e}}}\, \biggl[
{{\partial}\over{\partial\,\overline{x}}}\Bigl({{4}\over{3}}{{\partial\,
\overline{u}}\over{\partial\,\overline{x}}}
\,+\, {{1}\over{3}}{{\partial\,\overline{v}}\over{\partial\,\overline{y}}}
\Bigr) \,\,+\,\, {{1}\over{\epsilon^2}} \,\, 
{{\partial^{2}\,\overline{u}}\over{\partial\,\overline{y}^{2}}} \,\biggr] $}
\setbox30= \vbox {\halign{#\cr \box20 \cr \box21 \cr}}
\setbox31= \hbox{ $\vcenter {\box30} $}
\setbox44=\hbox{\noindent  (2.35) $\,\,   \left\{ \box31 \right. $}  
\smallskip \noindent $ \box44 $

\setbox20=\hbox{$\displaystyle 
{{\partial}\over{\partial\,\overline{t}}}(\overline{\rho}\,\overline{v})\,\,+\,\,
{{\partial}\over{\partial\,\overline{x}}}\bigl(\overline{\rho}\,\overline{u}
\,\overline{v}\,\bigr)\,\,+\,\,
{{\partial}\over{\partial\,\overline{y}}}\bigl(\overline{\rho}\,\overline{v}^{2}
\,+\, \,  {{1}\over{\epsilon^2}} \, \overline{p}\,\bigr)\,\,=$}
\setbox21=\hbox{$\displaystyle   \qquad \qquad  \qquad \qquad \qquad 
\,\,= \,\,  {{1}\over{\it R_{\rm e}}} \,
\biggl[\,{{\partial^{2}\,\overline{v}}\over{\partial\,\overline{x}^{2}}} \,\,+\,
\, \,  {{1}\over{\epsilon^2}} \,\, 
{{\partial}\over{\partial\,\overline{y}}}\Bigl({{1}\over{3}}{{\partial\,
\overline{u}}\over{\partial\,\overline{x}}}
\,+\,{{4}\over{3}}{{\partial\,\overline{v}}\over{\partial\,\overline{y}}} \Bigr)
\,\biggr]  $}
\setbox30= \vbox {\halign{#\cr \box20 \cr \box21 \cr}}
\setbox31= \hbox{ $\vcenter {\box30} $}
\setbox44=\hbox{\noindent  (2.36) $\,\,   \left\{ \box31 \right. $}  
\smallskip \noindent $ \box44 $

\setbox20=\hbox{$\displaystyle 
{{\partial}\over{\partial\,\overline{t}}} \biggl( \overline{\rho}\,\Bigl(
\overline{e}+{{1}\over{2}} \bigl( \overline{u}^{2}+ \,  \epsilon^{2} \overline{v}^{2}\bigr)
\Bigr) \biggr) \,+\,  $}
\setbox21=\hbox{$\displaystyle  \qquad   \,\,\, + \,\,\, \quad 
  {{\partial}\over{\partial\,\overline{x}}}\biggl(\overline{\rho}\,
\overline{u} \,\Bigl( \overline{e}+{{1}\over{2}} \bigl( \overline{u}^{2}+ \,  \epsilon^{2} 
\overline{v}^{2}\bigr) \Bigr)\,+\,\overline{p}\, \overline{u}\biggr)\,+$}
\setbox22=\hbox{$\displaystyle  \qquad   \,\,\, + \,\,\, \quad
  {{\partial}\over{\partial\,\overline{y}}}
\biggl(\overline{\rho}\,\overline{v} \,\Bigl( \overline{e}+{{1}\over{2}}  \bigl(
\overline{u}^{2}+ \,  \epsilon^{2}  \overline{v}^{2}\bigr) \Bigr)\,+\,
\overline{p}\,\overline{v}\biggr)\,\,=\,\, $}
\setbox23=\hbox{$\displaystyle \qquad \qquad  = \,\,\, \quad
 {{1}\over{\it R_{\rm e}}}\,\,
{{\partial}\over{\partial\,\overline{x}}} \,\biggl[\,\overline{u}\,\biggl(
{{4}\over{3}}{{\partial\,\overline{u}}\over{\partial\,\overline{x}}}
\,-\,{{2}\over{3}}{{\partial\,\overline{v}}\over{\partial\,\overline{y}}}
\biggr)\,+\,
\overline{v}\,\biggl({{\partial\,\overline{u}}\over{\partial\,\overline{y}}}
\,+ \,  \epsilon^{2} \, {{\partial\,\overline{v}}\over{\partial\,\overline{x}}}
\biggr) \biggl]   \,\,+  $}
\setbox24=\hbox{$\displaystyle  \qquad \qquad \quad  +\,\,\,
{{1}\over{\it R_{\rm e}}}\,\, {{\partial}\over{\partial\,\overline{y}}}
\,\biggl[\,\overline{u}\,\biggl( {{\partial\,\overline{v}}\over{\partial\,\overline{x}}}\,
+ \,  {{1}\over{\epsilon^2}} \,  \,
{{\partial\,\overline{u}}\over{\partial\,\overline{y}}} \biggr) \,+\, 
\overline{v}\,\biggl(
{{4}\over{3}}{{\partial\,\overline{v}}\over{\partial\,\overline{y}}}
\,-\,{{2}\over{3}}{{\partial\,\overline{u}}\over{\partial\,\overline{x}}}\biggr)\,
\biggl]  \,\,+ $}
\setbox25=\hbox{$\displaystyle \qquad \qquad  \quad +\,\,\,
{{1}\over{\it R_{\rm e}}}\,\,{{1}\over{\it P_{\rm r}}}\,\,\biggl(
{{\partial^{2}\,\overline{T}}\over{\partial\,\overline{x}^{2}}} \,+\,
{{1}\over{\epsilon^2}} \,  \,
{{\partial^{2}\,\overline{T}}\over{\partial\,\overline{y}^{2}}} \biggl)\,\,. $}
\setbox30= \vbox {\halign{#\cr  \box20 \cr \box21 \cr \box22 \cr \box23 \cr 
 \box24 \cr \box25 \cr}}
\setbox31= \hbox{ $\vcenter {\box30} $}
\setbox44=\hbox{\noindent  (2.37) $\,\,   \left\{ \box31 \right. $}  
\smallskip \noindent $ \box44 $

\bigskip \noindent $\bullet \quad$   
The boundary conditions associated with
these equations are of Dirichlet type on the boundary of the pipe~: 

\smallskip \noindent  (2.38) $\qquad \displaystyle
\overline{u}\,(x, y=-h) \,\,=\,\, \overline{u}\,(x, y=h) \,\,=\,\,0 $
\smallskip \noindent  (2.39) $\qquad \displaystyle
\overline{v}\,(x, y=-h)  \,\,=\,\,\overline{v}\,(x, y=h) \,\,=\,\,0 $
\smallskip \noindent  (2.40) $\qquad \displaystyle
\overline{T}\,(x, y=-h)  \,\,=\,\,\overline{T}\,(x, y=h)
\,\,=\,\,\overline{T}_{0}  $

\smallskip \noindent 
where $\,\overline{T}_{0}\,$ is the nondimensionless temperature given at the
boundary of the pipe.

\bigskip \noindent $\bullet \quad$
We observe first that the velocities  $\,\overline{u} \,$ and $\,\overline{v} \,$
have the same order of magnitude and in consequence, due to the fact that

\smallskip \noindent  (2.41) $\qquad \displaystyle
\epsilon^2 \,\,<<\,\,1 \hfill $

\smallskip \noindent
we can neglect in the left hand side of equation (2.37) the  $\,\overline{v} \,$
term compared to the $\,\overline{u} \,$ term. 

\bigskip \noindent $\bullet \quad$
We make the hypothesis that a typical distance for longitudinal variation of
physical fields is of the order $\,\Lambda\,$. In particular $\,\,\displaystyle
{{\partial\,u}\over{\partial\,x}}\,\simeq\, {{U}\over{\Lambda}}\,\,$  and in consequence

\smallskip \noindent  (2.42) $\qquad \displaystyle
{{\partial\,\overline{u}}\over{\partial\,\overline{x}}}\,\simeq\, \,1\,. $

\smallskip \noindent 
We observe that no bigger gradients than $\,\displaystyle {{1}\over{\Lambda}}\,$ are
taken in consideration into hypothesis (2.42) which means that the flow is regular
and that no turbulence occurs. In an analogous way, a typical distance for
transversal variation of all the fields is of the order of $\,l\,$ and in
particular  $\,\,\displaystyle {{\partial\,v}\over{\partial\,y}}\,\simeq\,
{{V}\over{l}}\,=\, {{U}\over{\Lambda}} \,$ and we have again

\smallskip \noindent  (2.43) $\qquad \displaystyle
{{\partial\,\overline{v}}\over{\partial\,\overline{y}}}\,\simeq\, \,1\,.$

\smallskip \noindent 
More generally all the differential expressions of the type   $\,\,\displaystyle 
{{\partial^{k}\,\overline{w}}\over{\partial\,\overline{z}^{k}}} \,$ with
$\,w\,$ equal to one of the physical fields $\,\rho, u, v, T, e\,$, variable
$\,z\,$ equal to $\, t, x, {\rm or} \, \, y\,$ and $\,k=1,2\,$, is finally of the
order of 1 (see {\it e.g.} Schlichting [Sc55] or Cousteix [Co88])~: 

\smallskip \noindent  (2.44) $\qquad \displaystyle
{{\partial^{k}\,\overline{w}}\over{\partial\,\overline{z}^{k}}}
\,\simeq\, \,1\,. $

\bigskip \noindent $\bullet \quad$ 
When we sum linear combinations of such expressions with coefficients of the type
$\,1\,$ or $\,  {{1}\over{\epsilon^2}} \,$ (as
in the right hand side of relation (2.35)), the leading term is the one that has
the dominant factor $\, {{1}\over{\epsilon^2}} \,$ as a coefficient. We
neglect in the following all the other terms. After  these approximations, we have derived
the so-called Thin Layer Navier Stokes equations  (see {\it e.g.}  Baldwin-Lomax [BL78], 
Kutler-Chakravarthy-Lombard [KCL78] or  Rubin and Tannehill [RT92]). We re-write them without any
adimensionalization~: 

\bigskip \noindent  (2.45) $\qquad \displaystyle
{{\partial\,\rho}\over{\partial\,t}}\,\,+\,\,{{\partial}\over{\partial\,x}}
(\rho\,u)\,\,+\,\,{{\partial}\over{\partial\,y}}(\rho\,v)\,\,=\,\,0 $

\smallskip \noindent  (2.46) $\quad \displaystyle
{{\partial}\over{\partial\,t}}(\rho\,u)\,\,+\,\,
{{\partial}\over{\partial\,x}}\bigl(\rho\,u^{2}\,+\,p\,\bigr)\,\,+\,\,
{{\partial}\over{\partial\,y}}\bigl(\rho\,u\,v\bigr)\,\,=\,\, \mu\, 
{{\partial^{2}\,u}\over{\partial\,y^{2}}} $

\smallskip \noindent  (2.47) $\quad \displaystyle
{{\partial p}\over{\partial\,y}}\,\,=\,\, \mu\,
{{\partial}\over{\partial\,y}}\,
\Bigl(\,{{1}\over{3}}{{\partial\,u}\over{\partial\,x}}
\,+\, {{4}\over{3}}{{\partial\,v}\over{\partial\,y}} \, \Bigr) $

\setbox21=\hbox{$\displaystyle 
{{\partial}\over{\partial\,t}} \biggl( \rho\,\Bigl(
e+{{1}\over{2}}  u^{2}\Bigr) \biggr)\,\,+\,\,
{{\partial}\over{\partial\,x}}\biggl(\rho\,u\,\Bigl(
e+{{1}\over{2}}  u^{2}\Bigr)\,+\,p\,u\biggr)\,\,+\,\,$}
\setbox22=\hbox{$\displaystyle   \qquad \qquad 
\,\,+ \,\,{{\partial}\over{\partial\,y}}\biggl(\rho\,v\,\Bigl(
e+{{1}\over{2}}u^{2} \Bigr)\,+\,p\,v\biggr)\,\,=\,\,
\mu\,{{\partial}\over{\partial\,y}}\Bigl(\,u\, {{\partial\,u}\over{\partial\,y}}\, 
\Bigr) + \,\, k\, {{\partial^{2}\,T}\over{\partial\,y^{2}}} \,.  $}
\setbox30= \vbox {\halign{#\cr \box21 \cr \box22 \cr}}
\setbox31= \hbox{ $\vcenter {\box30} $}
\setbox44=\hbox{\noindent  (2.48) $\,\,   \left\{ \box31 \right. $}  
\smallskip \noindent $ \box44 $

   \vfill\eject  
\noindent {\smcaps 3) $\quad$  Perfect fluid for main flow.} 
\smallskip \noindent $\bullet \quad$ 
We suppose now that the flow in the pipe satisfies the thin layer Navier-Stokes
equations (2.45)-(2.48) and for fixed time $\,t\,$ and abscissa $\,x\,$, we
integrate equation (2.45) between $\,y\,=\,-h\,$ and $\,y\,=\,+h\,$. We obtain in
this way

\setbox21=\hbox{$\displaystyle 
{{\partial}\over{\partial\,t}} \biggl( {{1}\over{2h}}
\int_{\displaystyle -h}^{\displaystyle {h}} \,\rho(t,x,y)\,{\rm d}y 
\biggr)\,\,+\,\, {{\partial}\over{\partial\,x}}\biggl( {{1}\over{2h}}
\int_{\displaystyle -h}^{\displaystyle {h}} \,(\rho u)(t,x,y)\,{\rm d}y  \biggr)\,\,
\,\,+\,\, $}
\setbox22=\hbox{$\displaystyle   \qquad \qquad  \qquad \qquad \qquad 
\,\,+  \,\, {{1}\over{2h}} \, \Bigl( (\rho v)(t,x,h) \,-\, (\rho v)(t,x,-h) 
\Bigr)\,\,=\,\,0\,.  $}
\setbox30= \vbox {\halign{#\cr \box21 \cr \box22 \cr}}
\setbox31= \hbox{ $\vcenter {\box30} $}
\setbox44=\hbox{\noindent  (3.1) $\,\,   \left\{ \box31 \right. $}  
\smallskip \noindent $ \box44 $

\smallskip \noindent
Due to boundary condition (2.39), the third term in (3.1) is null. We introduce
now the mean values of density, momentum and energy in each $x\,$ section
according to~:

\smallskip \noindent  (3.2) $\qquad \displaystyle
\widetilde{\rho} (t,x)\,\,=\,\,{{1}\over{2h}}
\int_{\displaystyle -h}^{\displaystyle {h}} \,\rho(t,x,y)\,{\rm d}y  $

\smallskip \noindent  (3.3) $\qquad \displaystyle
\widetilde{\rho} (t,x)\,\widetilde{u} (t,x)  \, \,=\,\,{{1}\over{2h}}
\int_{\displaystyle -h}^{\displaystyle {h}} \,(\rho u)(t,x,y)\,{\rm d}y  $

\smallskip \noindent  (3.4) $\qquad \displaystyle 
\widetilde{\rho} (t,x) \, \widetilde{e} (t,x) \, +\, {{1}\over{2}} \widetilde{\rho}
(t,x)\,\widetilde{u}^{2} (t,x)  \, \,=\,\,{{1}\over{2h}}
\int_{\displaystyle -h}^{\displaystyle {h}} \,\bigl( \, \rho e\,+\,{{1}\over{2}}
\rho u^{2} \, \bigr) (t,x,y)\,{\rm d}y  \,.$

\smallskip \noindent 
In terms of these new variables, the conservation of mass stands as~: 

\smallskip \noindent  (3.5) $\qquad \displaystyle 
{{\partial\,\widetilde{\rho}}\over{\partial t}}\,\,+\,\,
{{\partial}\over{\partial x}} \bigl(\, \widetilde{\rho} \, \widetilde{u}  \,
\bigr)\,\,=\,\,0 \,.$

\bigskip \noindent $\bullet \quad$ 
In a similar way we integrate the impulse and the energy equations in the thickness
of the pipe. We get 

\setbox21=\hbox{$\displaystyle   
{{\partial}\over{\partial\,t}}\bigl(\, \widetilde{\rho} \, \widetilde{u}  \,
\bigr)\,\,+\,\, {{\partial}\over{\partial\,x}}  \biggl( {{1}\over{2h}}
\int_{\displaystyle -h}^{\displaystyle h} \, \bigl( \, \rho
u^{2}\,+\,p\, \bigr) (t,x,y) \,{\rm d}y  \biggr)  \,\,=$}
\setbox22=\hbox{$\displaystyle   \qquad \qquad  \qquad \qquad  \qquad
\qquad  \qquad \,\,=\,\, {{\mu}\over{2\,h}}  \,\, \biggl(\,{{\partial u}\over{\partial
y}}\,(t,x,h)\,-\,  {{\partial u}\over{\partial y}}\,(t,x,-h) \, \biggr) $}
\setbox30= \vbox {\halign{#\cr \box21 \cr \box22 \cr}}
\setbox31= \hbox{ $\vcenter {\box30} $}
\setbox44=\hbox{\noindent  (3.6) $\,\,   \left\{ \box31 \right. $}  
\smallskip \noindent $ \box44 $

\setbox21=\hbox{$\displaystyle   
{{\partial}\over{\partial\,t}}\bigl(\, 
\widetilde{\rho} \, \widetilde{e}  \, +\, {{1}\over{2}} \widetilde{\rho} \, 
\widetilde{u}^{2} \,\bigr)\,+\,{{\partial}\over{\partial\,x}}
\biggl( {{1}\over{2h}} \int_{\displaystyle -h}^{\displaystyle h} \, 
\Bigl( \, \rho  u \, \bigl( \, e+{{u^{2}}\over{2}} \, \bigr) +pu\,
\Bigr) (t,x,y) \,{\rm d}y  \biggr)  \,=     $}
\setbox22=\hbox{$\displaystyle   \qquad \qquad  \qquad \qquad  \qquad
\qquad  \qquad \,\,=\,\,  {{k}\over{2\,h}}  \,\,\biggl(\,{{\partial T}\over{\partial
y}}\,(t,x,h)\,-\,  {{\partial T}\over{\partial y}}\,(t,x,-h) \, \biggr)   $}
\setbox30= \vbox {\halign{#\cr \box21 \cr \box22 \cr}}
\setbox31= \hbox{ $\vcenter {\box30} $}
\setbox44=\hbox{\noindent  (3.7) $\,\,   \left\{ \box31 \right. $}  
\smallskip \noindent $ \box44 $

\smallskip \noindent
due to boundary conditions (2.38) and (2.39).

\bigskip \noindent $\bullet \quad$
We make first the hypothesis that the fields are quasi-constant in each section of the
pipe and second that they have a rapid variation in a boundary layer region of thickness
$\,\delta
\,$, with the condition 

\smallskip \noindent  (3.8) $\qquad \displaystyle 
S_{h} \,\,\equiv \,\,{{\delta}\over{h}} \,\, \, << 1 \,. $

\smallskip \noindent 
The first  hypothesis is absolutly non trivial. 

\bigskip \noindent $\bullet \quad$
In aerodynamics, it conducts (see {\it e.g.} Whitham [Wh74] and Msallam [Ms98])  to the shallow
water equations when the following physical hypothesis is satisfied~: 

\smallskip \noindent  (3.9) $\qquad \displaystyle 
{{1}\over{\it R_{\rm e}^{\rm aero}}} \,{{1}\over{S_{h}}}\,\,\approx\,\,{{\delta}\over{L}}
\, {{h}\over{L}}  \,\, = \,\, \epsilon \,{{h}\over{L}}  \,\,\, << \,1 $

\smallskip \noindent 
due to hypothesis (2.1) and choice of variable  $\,\epsilon\,$ done in (2.5).

\bigskip \noindent $\bullet \quad$
In acoustics, the Reynolds number is modified according to relation (2.9),  {\it i.e.}

\smallskip \noindent  (3.10) $\qquad \displaystyle 
{\it R}^{\rm acou}_{\rm e}\,\,= \,\, {{\rho_{0}\,c_{0}\,\lambda}\over{2\,\pi\,\mu}}$

\smallskip \noindent  
and the hypothesis of quasi-constancy of all the fields in the main flow is
satisfied under the hypothesis (Kergomard [Ke81], Menguy and Gilbert [MG97])

\smallskip \noindent  (3.11) $\qquad \displaystyle 
{{1}\over{\it R}^{\rm acou}_{\rm e}} \, {{1}\over{S_{h}}} \,\, << \,1 \,. $

\bigskip \noindent $\bullet \quad$ 
This hypothesis can be justied as follow.
Observe first that for a simple linear wave, the variation of pressure
$\,p_{a}\,$ due to acoustics  satisfies the relation 

\smallskip \noindent  (3.12) $\qquad \displaystyle 
p_{a} \,\,=\,\, \rho_{0}\,c_{0}\,u\,\,. $

\smallskip \noindent 
Second the transverse gradient of pressure $\displaystyle \,{{\partial p }\over
{\partial y}}\,$ satisfies at the first order a linearized version of equation
(2.36)~: 

\setbox20=\hbox{$\displaystyle  \rho_{0}\, {{\partial v }\over {\partial t}} $}
\setbox27=\hbox{$\displaystyle  {{\partial p }\over {\partial y}} $}
\smallskip \noindent  (3.13) $\qquad \displaystyle 
\mod{\box27} \,\,=\,\, \mod{\box20} $

\smallskip \noindent 
and the right hand side of this equation (3.13) can be evaluated as follow~: 

\smallskip \noindent  (3.14) $\qquad \displaystyle 
{{\partial p }\over {\partial y}} \quad \approx \quad \rho_{0} \,\,\, {{\delta}
\over {\lambda}} \, u \,\,\, {{1}\over{\tau}}  \quad = \quad 
\rho_{0} \, u \,\,  c_{0} \, \Bigl( {{\delta} \over {\lambda}} \Bigr) ^{2} \,\,
{{1}\over{\delta}}  $

\smallskip \noindent 
because $\, \displaystyle \tau = {\lambda \over {c_0}} \,$ 
where the thickness of the boundary layer $\, \delta \,$ is (see {\it e.g.}
Lighthill [Li78]) of the order of 
$\displaystyle \, {{\lambda}\over{\sqrt{{\it R}^{\rm acou}_{\rm e}}}} \,$~: 

\smallskip \noindent  (3.15) $\qquad \displaystyle
\overline{\delta} \,\, = \,\, {{\delta} \over {\lambda}} \,\,
\approx \,\,{{1}\over { \sqrt {{\it R}^{\rm acou}_{\rm e}}}} \,. $

\smallskip \noindent 
We insert relations (3.8), (3.12) and (3.15) inside (3.14) and obtain 

\smallskip \noindent  (3.16) $\qquad \displaystyle
{{\partial p }\over {\partial y}} \,\, \approx \,\, {{p_{a}}\over{h}} \, 
{{1}\over{{\it R}^{\rm acou}_{\rm e}}} \, {{1}\over{S_{h}}} \,. $

\smallskip \noindent 
The transverse variations of pressure are of the order of the axial variations of
pressure multiplied by the factor $\displaystyle \, {{1}\over{{\it R}^{\rm acou}_{\rm e}}
\, {S_{h}}}$ . Then relation (3.11) express that the transverse
variation of pressure can be neglected compared with the axial ones. 

\bigskip \noindent $\bullet \quad$
We observe also that 

\smallskip \noindent  (3.17) $\qquad \displaystyle
{{1}\over{\it R}^{\rm acou}_{\rm e}} \, {{1}\over{S_{h}}} \,\,=\,\,
\Bigl( {{\delta}\over{\lambda}} \Bigr)^2 \, {{h}\over{\delta}} \,\,=\,\,
{{\delta \, h}\over{\lambda^2}} \,\,= \,\, 
{{\delta}\over{\lambda}} \, \, {{h}\over{\lambda}}  \,\, \leq \,\, 
{1\over4} \, {{1}\over { \sqrt {{\it R}^{\rm acou}_{\rm e}}}}  \,\, << \,\, 1 \, $

\smallskip \noindent 
due to hypotheses (2.6), (3.8) and (3.15). Then   hypothesis (3.11) is established even  if
the boundary layer thickness $\,\delta\,$ is greater than the order of magnitude of the
characteristic length $\,l\,$  of visco-thermic effects. Physically, it corresponds to
neglect volume losses compared to wall losses.

\bigskip \noindent $\bullet \quad$ 
Under the hypothesis that all the fields are
constant in the section  

\smallskip \noindent   $\qquad \qquad \displaystyle
-(h-\delta)\,\,\leq \,\,y \,\,\leq \,\, (h-\delta)\,, $

\smallskip \noindent
we first observe that pressure $\,\widetilde{p}\,$ associated via the state equation
(2.13) to mean density $\,\widetilde{\rho}\,$ and mean internal energy 
$\,\widetilde{e}\,$ can be well approached by the mean value of pressure~: 

\smallskip \noindent  (3.18) $\qquad \displaystyle
(\gamma -1)\,\widetilde{\rho}\,\widetilde{e}
\,\,\,\simeq\,\,\, {{1}\over{2h}} \int_{\displaystyle -h}^{\displaystyle h} 
\,\,p(t,x,y) \,{\rm d}y  $

\smallskip \noindent 
as mentioned previously. In an analogous way, we have

\smallskip \noindent  (3.19) $\qquad \displaystyle
\widetilde{\rho}\,\widetilde{u}^{2}
\,\,\,\simeq\,\,\, {{1}\over{2h}} \int_{\displaystyle -h}^{\displaystyle h} 
\bigl( \rho \, u^{2} \bigr)\,(t,x,y) \,{\rm d}y  $
\smallskip \noindent  (3.20) $\qquad \displaystyle
\widetilde{\rho} \, \widetilde{u}\, \bigl( \widetilde{e}  \, +\, {{1}\over{2}}
\widetilde{u}^{2} \bigr) \,+\,\widetilde{p}\,\widetilde{u}
\,\,\,\simeq\,\,\, {{1}\over{2h}} \int_{\displaystyle -h}^{\displaystyle h} 
\Bigl(  \rho  u \, \bigl( \, e+{{u^{2}}\over{2}} \, \bigr) +pu\,
\Bigr) (t,x,y) \,{\rm d}y \,.$

\smallskip \noindent 
All the hypotheses (3.18)-(3.20) suppose finally that mean values of a nonlinear
 function is quasi-equal to the same nonlinear function of the mean values. This
hypothesis is correct when the nonlinear function is well approximated by a
constant. 

\bigskip \noindent $\bullet \quad$
We can now insert relations (3.18) to (3.20) inside equations (3.6) and (3.7). We
obtain the final model for unidimensional perfect flow~: 

\smallskip \noindent  (3.21) $\quad \displaystyle
\widetilde{p}(t,x)\,\,\,\equiv \,\,\,(\gamma -1)\,\widetilde{\rho}\,\widetilde{e}$
\smallskip \noindent  (3.22) $\quad \displaystyle
{{\partial\,\widetilde{\rho}}\over{\partial t}}\,\,+\,\,
{{\partial}\over{\partial x}} \bigl(\, \widetilde{\rho} \, \widetilde{u}  \,
\bigr)\,\,=\,\,0 \,.$
\smallskip \noindent  (3.23) $\quad \displaystyle
{{\partial}\over{\partial\,t}}\bigl(\, \widetilde{\rho} \, \widetilde{u}  \,
\bigr)\,\,+\,\, {{\partial}\over{\partial\,x}}  \biggl(
\widetilde{\rho}\,\widetilde{u}^{2}\,+\, \widetilde{p} \biggr) \,\,=
 {{\mu}\over{2\,h}} \,\, \biggl(\,{{\partial u}\over{\partial y}}\,(t,x,h)\,-\, 
{{\partial u}\over{\partial y}}\,(t,x,-h) \, \biggr) $

\setbox21=\hbox{$\displaystyle   
{{\partial}\over{\partial\,t}}\bigl(\, 
\widetilde{\rho} \, \widetilde{e}  \, +\, {{1}\over{2}} \widetilde{\rho} \, 
\widetilde{u}^{2} \,\bigr)\,+\,{{\partial}\over{\partial\,x}} \biggl(
\widetilde{\rho} \, \widetilde{u}\, \bigl( \widetilde{e}  \, +\,
{{1}\over{2}} \widetilde{u}^{2} \bigr) \,+\,\widetilde{p}\,\widetilde{u}
\biggr)  \,=     $}
\setbox22=\hbox{$\displaystyle   \qquad \qquad  \qquad \qquad  
\qquad  \qquad  \,\,=\,\,  {{k}\over{2\,h}} \,\,\biggl(\,{{\partial T}\over{\partial
y}}\,(t,x,h)\,-\,  {{\partial T}\over{\partial y}}\,(t,x,-h) \, \biggr)  \,\,. $}
\setbox30= \vbox {\halign{#\cr \box21 \cr \box22 \cr}}
\setbox31= \hbox{ $\vcenter {\box30} $}
\setbox44=\hbox{\noindent  (3.24) $\,\,   \left\{ \box31 \right. $}  
\smallskip \noindent $ \box44 $

\bigskip
\bigskip
  \noindent {\smcaps 4) $\quad$  Acoustic boundary layer.}
\smallskip \noindent $\bullet \quad$ 
We suppose as previously that the flow in the pipe satisfies the Thin Layer Navier
Stokes equations (2.45)-(2.48). In the following, we look for  the boundary
layer $(\, y \leq -h + \delta \, $ or $ \, y \geq h-\delta). $ The previous Thin
Layer Navier Stokes equations   simplify and we  obtain the equations of acoustic
boundary layer [Ch64].

\bigskip \noindent $\bullet \quad$ 
First we suppose that some reference state with null velocity is given. It is
a priori the air at usual atmospheric pressure $\,p_0\,$ and usual temperature
$\,\theta_0$ ;   with the state low of perfect gas, the reference density
$\,\rho_0\,$ is  given.  Second we search a field $\,
\rho(t,\,x,\,y),\,  u(t,\,x,\,y),\,  T(t,\,x,\,y),\,  p(t,\,x,\,y) $ of the form

\setbox20=\hbox{$\displaystyle 
\rho (t,\,x,\,y) \,\,=\,\, \rho_0 \,+\, \rho' (t,\,x,\,y)$}
\setbox21=\hbox{$\displaystyle  
u (t,\,x,\,y)\,\,=\,\, 0 \, \, \, +\, u' (t,\,x,\,y) $}
\setbox22=\hbox{$\displaystyle 
T (t,\,x,\,y) \,\,=\,\, \theta_0 \,+\, T' (t,\,x,\,y) $}
\setbox23=\hbox{$\displaystyle 
p (t,\,x,\,y) \,\, = \,\, p_0 \, \, +\, p' (t,\,x,\,y) \,.$}
\setbox30= \vbox {\halign{#\cr \box20 \cr \box21 \cr \box22 \cr \box23 \cr}}
\setbox31= \hbox{ $\vcenter {\box30} $}
\setbox44=\hbox{\noindent  (4.1) $\,\,   \left\{ \box31 \right. $}  
\smallskip \noindent $ \box44 $

\smallskip \noindent
We linearize the equations (2.45)-(2.48) around the reference state $\,( \rho_0, \,
p_0, \,\theta_0).$ We suppose that we can neglect the nonlinear contributions of
the bondary layer for studing the equations (3.21)-(3.24) of the main flow. We
recall that these equations have been obtained by integrating the Thin Layer Navier
Stokes equations on the complete width of the channel and the detailed analysis of
the different contributions have been derived by Msallam [Ms98]. Recall that by
doing this linear approximation, we neglect the acoustic streaming effect (see {\it e.g.}
Batchelor [Ba67], Makarov and Ochmann [Mo97]), all separation effects inside the
boundary layer (Merkli and Thoman [MT75]) and all random unstationary effects of
turbulence [MT75]. The algebra is classical (Chester [Ch64]) and straightforward.
We obtain  

\smallskip \noindent  (4.2) $\qquad \displaystyle
{{p'}\over{p_0}} \,\,=\,\, {{\rho'}\over{\rho_0}} \,+\, {{T'}\over{\theta_0}}$
\smallskip \noindent  (4.3) $\qquad \displaystyle
{{\partial \rho'}\over{\partial t}} \,\,+\,\, \rho_0 \, \, {\rm div}\,u'
\,\,=\,\,0.$
\smallskip \noindent  (4.4) $\qquad \displaystyle
\rho_0 \, {{\partial u'}\over{\partial t}} \,-\, \mu \,{{\partial^2 u'}
\over{\partial y^2}} \,\,=\,\,-\, {{\partial p'}\over{\partial x}} $
\smallskip \noindent  (4.5) $\qquad \displaystyle
\rho_0 \, C_p \, {{\partial T'}\over{\partial t}} \,-\, k \,{{\partial^2 T'}
\over{\partial y^2}} \,\,=\,\, {{\partial p'}\over{\partial t}}  \,.$

\smallskip \noindent 
We note (with Chester [Ch64]) that equations (4.4) and (4.5) are two heat
equations coupled via the right hand sides. 

\bigskip \noindent $\bullet \quad$ 
The boundary conditions  associated to (say) the bottom of the boundary layer
($y=-h$) are~:  

\smallskip \noindent  (4.6) $\qquad \displaystyle
u'(t,\,x,\,y=-h) \,\,= \,\, 0.$
\smallskip \noindent  (4.7) $\qquad \displaystyle
T'(t,\,x,\,y=-h) \,\,= \,\, 0.$

\smallskip \noindent 
At the top of the boundary layer ($\, y \approx -h + \delta \,$), we must mutch
the boundary layer flow with the main flow~:  

\smallskip \noindent  (4.8) {$\qquad \displaystyle
u'(t,\,x,\,y \longrightarrow -h+\delta) \,\,\longrightarrow $ (velocity in the main
flow)$ (t,\,x) $ } 
\smallskip \noindent  (4.9) {$\qquad \displaystyle
T'(t,\,x,\,y \longrightarrow -h+\delta) \,+\, \theta_0 \,\,\longrightarrow $
(temperature  in the main flow)$ (t,\,x) \,. $ }

\bigskip
\bigskip
 \noindent {\smcaps 5) $\quad$  The coupled problem.}
\smallskip \noindent $\bullet \quad$ 
We couple in this section the main flow in the pipe described in section 3 with
the acoustic boundary layer presented in part 4. More precisely, the main flow is
described by three unknown functions (density, velocity, internal energy)~: 

\smallskip \noindent  (5.1) $\,  \displaystyle
[0, \, +\infty [  \times  [0, \, L ] \ni (t,\,x) \longmapsto
\bigl( \rho(t,x), \, u(t,x), \, e(t,x) \bigr) \in [0, \, +\infty [  \times 
\R  \times  [0, \, +\infty [  $

\smallskip \noindent
which represent the mean value in the section of the pipe of each field (denoted
with a tilda in section 3). In the bounday layer, we suppose that the faces $\,y =
\pm h \,$ are composed by symmetric flows and we fix some transverse variable
$\,\eta \in [0, \, +\infty [$. 

\bigskip 
\centerline { \epsfysize=3cm    \epsfbox  {fig2.epsf} }
\smallskip  \smallskip
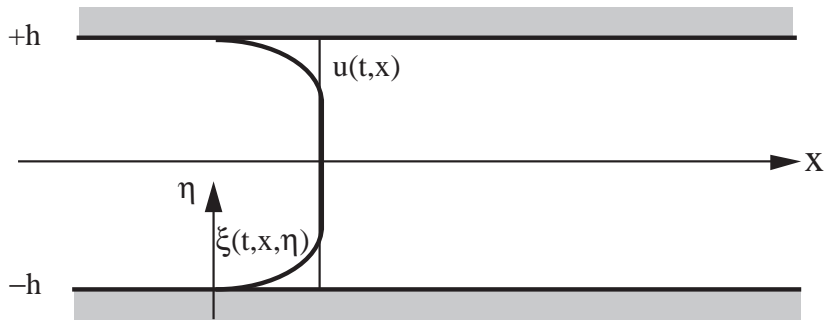
\centerline {{\bf Figure 2}	\quad Velocity field $\, u(x,\,t) \,$ in the mean flow}
\centerline { and  velocity field $\, \xi(t,\,x,\,\eta)\,$ inside the boundary layer. }
\smallskip

\smallskip \noindent
The unknowns are velocity $\, \xi \,$ and
temperature $\, \theta \,$ in the boundary layer~:  

\smallskip \noindent  (5.2) $  \displaystyle
[0, \, +\infty [  \times  [0, \, L ] \times  [0, \, +\infty [ 
\ni \! (t,\,x, \, \eta ) \longmapsto 
\bigl( \xi (t,\,x, \, \eta ), \, \theta  (t,\,x, \, \eta ) \bigr) \!
\in  \R \times [0, \, +\infty [  $

\smallskip \noindent
Notice the important point concerning the modelling~:  we consider on one hand {\bf
two} velocity fields $\, u  \,$ and $\, \xi \,$ and on the other hand {\bf two}
temperature fields $\displaystyle \, T \, = \, e/C_v\,$ and $\, \theta$.

\bigskip \noindent $\bullet \quad$ 
The transverse scale $\, \eta \,$ for describing the boundary layer flow is very
small compared to the transverse dimension $\,h\,$ of the flow. Then it is
consistent to set boundary conditions for $\, \eta \, \longrightarrow \, +\infty 
\,$~:

\setbox20=\hbox{$\displaystyle 
{{\partial \xi}\over{\partial \eta}}(t,\,x, \, \eta )  \, \longrightarrow \,
0 \qquad $when$ \,\, \eta \, \longrightarrow \, +\infty  $}
\setbox21=\hbox{$\displaystyle  
{{\partial \theta}\over{\partial \eta}} (t,\,x, \,
\eta )  \, \longrightarrow \,0 \qquad $when$ \,\, \eta \, \longrightarrow \,
+\infty \, $}
\setbox30= \vbox {\halign{#\cr \box20 \cr \box21 \cr}}
\setbox31= \hbox{ $\vcenter {\box30} $}
\setbox44=\hbox{\noindent  (5.3) $\quad   \left\{ \box31 \right. $}  
\smallskip \noindent $ \box44 $

\smallskip \noindent
and we have done this particular choice in our simulations. Nevertheless, stronger boundary
conditions for  $\, \eta \, \longrightarrow \, +\infty  \,$    that are compatible with
the  observed solutions in our numerical experiments could be the following ones~:

\setbox20=\hbox{$\displaystyle 
\xi (t,\,x, \, \eta )  \, \longrightarrow \, u(t,\,x) \qquad $when$ \,\,
\eta \, \longrightarrow \, +\infty  $}
\setbox21=\hbox{$\displaystyle  
\theta (t,\,x, \, \eta )  \, \longrightarrow \, T(t,\,x) \qquad $when$ \,\,
\eta \, \longrightarrow \, +\infty  \,.$}
\setbox30= \vbox {\halign{#\cr \box20 \cr \box21 \cr}}
\setbox31= \hbox{ $\vcenter {\box30} $}
\setbox44=\hbox{\noindent  (5.4) $\quad   \left\{ \box31 \right. $}  
\smallskip \noindent $ \box44 $

\bigskip \noindent $\bullet \quad$ 
For $\, \eta \,= \, 0$, we just have to consider Dirichlet boundary conditions 

\setbox20=\hbox{$\displaystyle 
\xi(t,\,x,\,0) \,\,=\,\,0 $}
\setbox21=\hbox{$\displaystyle  
\theta(t,\,x,\,0) \,\,=\,\,\theta_0$}
\setbox30= \vbox {\halign{#\cr \box20 \cr \box21 \cr}}
\setbox31= \hbox{ $\vcenter {\box30} $}
\setbox44=\hbox{\noindent  (5.5) $\quad   \left\{ \box31 \right. $}  
\smallskip \noindent $ \box44 $

\smallskip \noindent
where $\, \theta_0 \,$ is the value of imposed temperature on the walls. 

\bigskip \noindent $\bullet \quad$ 
The partial differential equations for the evolution of main flow variables are
simply derived from equations (3.21)-(3.24) ; the source terms of stress  and
thermal flux  at the wall in the right hand side of equations (3.23) and (3.24) are 
nomore obtained by solving the Thin Layer Navier Stokes equations
(2.45)-(2.48) but the ones coming from the boundary layer model (4.2)-(4.5).
We obtain in this way~:  

\smallskip \noindent  (5.6) $\qquad \displaystyle
p(t,x)\,\,\,\equiv \,\,\,(\gamma -1)\,\rho\,e$
\smallskip \noindent  (5.7) $\qquad \displaystyle
{{\partial\,\rho}\over{\partial t}}\,\,+\,\,
{{\partial}\over{\partial x}} \bigl(\, \rho \, u  \, \bigr)\,\,=\,\,0 \,$
\smallskip \noindent  (5.8) $\qquad \displaystyle
{{\partial}\over{\partial\,t}}\bigl(\, \rho \, u  \, \bigr)\,\,+\,\,
{{\partial}\over{\partial\,x}}  \bigl(\, \rho\,u^{2}\,+\, p \, \bigr) \,\,= \,\, 
-  {{\mu}\over{h}} \,\,{{\partial \xi}\over{\partial \eta}}\,\bigl(t,\,x,\, 0\bigr)
\, $
\smallskip \noindent  (5.9) $\qquad \displaystyle
{{\partial}\over{\partial\,t}} \bigl(\, \rho \, e  \, 
 \, +\,\,  {{1}\over{2}}\, \rho \, u^2 \,  \bigr)\,\,+\,\,
{{\partial}\over{\partial\,x}}  \bigl( \, \rho \, u \, e  \, 
 \, +\,\, {{1}\over{2}}\, \rho \, u^3 \,   \, +\,  \, p \, u \, \bigr)
\,\,= \,\,  -  {{k}\over{h}} \,\,{{\partial \theta}\over{\partial
\eta}}\,\bigl(t,\,x,\,0 \bigr) \,. $

\bigskip \noindent $\bullet \quad$ 
The evolution equations for boundary layer variables are obtained in a similar way
from the heat equations (4.4)-(4.5) that model an acoustic boundary layer. With our
coupled model, the pressure term comes from the main one-dimensional model and no
more from the Thin  Layer Navier Stokes equations ; it is therefore considered as a
source term and for this reason is placed on the right hand side of the equations.
We get 

\smallskip \noindent  (5.10) $\qquad \displaystyle
\rho_0 \,  {{\partial \xi}\over{\partial t}} \,\,-\,\, \mu
{{\partial^2 \xi}\over{\partial \eta^2}} \, \,\,=\,\, -\,{{\partial p} \over
{\partial x}} $
\smallskip \noindent  (5.11) $\qquad \displaystyle
\rho_0 \,C_p \,  {{\partial \theta}\over{\partial t}} \,\,-\, \, k 
{{\partial^2 \theta}\over{\partial \eta^2}} \, \,\,=\,\, {{\partial p}
\over {\partial t}} \,.$

\bigskip \noindent $\bullet \quad$ 
We observe that equation (5.10) is a dynamic equation that allow the prediction
of velocity field $\, \xi(t,\,x,\,\eta)\,$ as long as pressure field
$\,p(t,\,x)\,$ is known. It is not so clear for temperature equation (5.11) for
the variable $\, \theta\,$ due to the dymamic term $  \,{{\partial p}
\over {\partial t}} \,$ on the right hand side. Nevertheless, following a
remark of Brenier [Br97], system (5.5) (second equation) and (5.11) can be replaced
by the new unknow function $\,
\sigma \,$ 

\smallskip \noindent  (5.12) $\qquad \displaystyle
\sigma \,\, = \,\, \rho_0 \, Cp \, \theta \,-\, p$

\smallskip \noindent
that satisfies the following heat equation 

\smallskip \noindent  (5.13) $\qquad \displaystyle
{{\partial \sigma}\over{\partial t}} \,\,-\, \, {{k} \over {\rho_0 \, Cp}} \, 
{{\partial^2 \sigma}\over{\partial \eta^2}} \, \,\,=\,\,0\,$

\smallskip \noindent  
due to the fact that   $\,\displaystyle {{\partial p} \over {\partial \eta}} \, \equiv
\, 0 .$  The following nonhomogeneous Dirichlet boundary condition is valid at the bottom
of the boundary layer~:  

\smallskip \noindent  (5.14) $\qquad \displaystyle
\sigma (t,\,x,\,0) \,\,=\,\,\,\, \rho_0 \, Cp \, \theta_0 \,-\, p(t,x) \,.$

\bigskip \noindent $\bullet \quad$ 
It is clear that on one hand, that stress viscous term $ \, \mu \,\smash{
{{\partial \xi} \over {\partial \eta}}}(\eta=0)\,$ and thermal flux $ 
\,k \, \smash{{{\partial \theta} \over {\partial \eta}}}(\eta=0)\,$ forces the main flow
equations (5.7)-(5.9) and on the other hand that  pressure field in the inviscid
flow forces the boundary layer equations (5.10)-(5.11).  We insist again on the
fact that the main originality of our  coupled model (5.6)-(5.11) consists of
choosing {\bf two} independent  unknowns functions for velocity (in the main flow
and in the bounday layer) and also {\bf two} independent  unknowns functions for
temperature. We do not make  the tentative to determine explicitely the boundary
layer thickness $\,\delta(t,x) \,$ or the displacement thickness  $\,\delta^*(t,x) \,$
(see {\it e.g.} Le Balleur [LB80]) in the way we set the coupled problem. In our approach,
the  boundary layer  thickness for momentum and energy  can be evaluated as a
global (and nontrivial) result from the entire knowledge of functions $\,\,  [0,
\, +\infty [  \, \ni \eta  \longmapsto  \xi (t,\,x, \, \eta ) \,$ and $\, [0, \,
+\infty [  \, \ni \eta  \longmapsto  \theta (t,\,x, \,\eta )$.

\bigskip \noindent $\bullet \quad$ 
We recall briefly also the inflow-outflow boundary conditions at $\, x = 0 \,$ and
$\, x = L \,$  concerning  the mean flow variables. At the inflow
($x=0$), the flow is subsonic then two conditions have to be considered~:  for axample, we
give on one hand some data concerning the input velocity field $\,u_0(t)\,$ or the input 
pressure field $\, \pi_0(t) \,$ and on the other hand the fact that entropy is not
dissipated at the entrance of the channel~:  

\smallskip \noindent  (5.15) $\qquad \displaystyle
u(t \,,\, 0) \,\,=\,\, u_0(t) \qquad $or$ \qquad p(t \,,\, 0) \,\,=\,\,\pi_0(t) $ 
\smallskip \noindent  (5.16) $\qquad \displaystyle
{{\partial}\over{\partial t}} \Bigl( {{p}\over{\rho^{\gamma}}} \Bigr) 
(t \,,\, 0) \,\,=\,\, 0 \,. $

\smallskip \noindent
At the outflow, the mean field remains subsonic and the theory of characteristics
(see {\it e.g.} Kreiss [Kr70]) show that  only one scalar boundary condition is
sufficient to set correctly the problem ; we choice nonreflecting boundary
conditions (see Whitham [Wh74] or Hedstrom [He79])~:   the outgoing wave is a
so-called $\,C_+\,$  simple wave [Wh74] {\it i.e.} both specific entropy $\,S \equiv
\displaystyle  {{p}\over{\rho^{\gamma}}} \,$ and Riemann invariant $\,\displaystyle R
\equiv  u-{{2c}\over{\gamma-1}} \,$ take constant values everywhere in this part of the
flow. In consequence this Riemann invariant $\,R\,$ is  advected  with all characteristic
celerities
 without distorsion and in particular the one with $\,u-c\,$ velocity  ~: 

\smallskip \noindent  (5.17) $\qquad \displaystyle
\biggl( {{\partial}\over{\partial t}} \Bigl( u-{{2c}\over{\gamma-1}} \Bigr) \,
\,+\,\, (u-c) \, {{\partial}\over{\partial x}} \Bigl( u-{{2c}\over{\gamma-1}}
\Bigr)  \biggr) \, (t \,,\, L) \,\,=\,\, 0 \,.$

\bigskip
\bigskip
 \noindent {\smcaps 6) $\quad$  Generalization to axisymmetric geometry.}
\smallskip \noindent $\bullet \quad$ 
In this section, the pipe is nomore a two-dimensional channel but a three-dimensional
cylinder  with an axisymmetric geometry. The length of the pipe is still denoted by $\, L
\,$ and the letter $\,h \, $ is used for the radius instead of half of the section.
Hypothesis (2.1) concerning ratio $\, \displaystyle {{h}\over{L}} \,$
remains valid and we have 

\smallskip \noindent  (6.1) $\qquad \displaystyle
{{h}\over{L}} \quad << \quad 1 \, . \,$

\smallskip \noindent
As is section 2, Thin Layer Navier Stokes equations are a good approximation of the flow
inside the entire geometry and  this model takes  now the following algebraic form in
this axisymmetric geometry~:   

\smallskip \noindent  (6.2) $\qquad \displaystyle
{{\partial\,\rho}\over{\partial\,t}}\,\,+\,\,{{\partial}\over{\partial\,x}}
(\rho\,u)\,\,+\,\,{{1}\over{y}}\,  {{\partial}\over{\partial\,y}}
\bigl( \rho \, v \,y \bigr) \,\,=\,\,0 $

\smallskip \noindent  (6.3) $\qquad \displaystyle
{{\partial}\over{\partial\,t}}(\rho\,u)\,\,+\,\,
{{\partial}\over{\partial\,x}}\bigl(\rho\,u^{2}\,+\,p\,\bigr)\,\,+\,\,
{{1}\over{y}}\, {{\partial}\over{\partial\,y}}\bigl(\rho\,u\,v \, y \bigr)\,\,=\,\, 
{{\mu}\over{y}}\, {{\partial}\over{\partial\,y}} \, \biggl( y \, 
{{\partial u}\over{\partial\,y}} \biggr) \,  $

\smallskip \noindent  (6.4) $\qquad \displaystyle
{{1}\over{y}}\, {{\partial}\over{\partial\,y}}\bigl(\, \rho\,v^{2}\,y\,+\,p\,y\, \bigr)\,\,=
 \,\, {{\mu}\over{y}}\, {{\partial}\over{\partial\,y}} \, \biggl(
\,{{1}\over{3}}{{\partial\,u}\over{\partial\,x}}
\,+\, {{4 \, y}\over{3 }}{{\partial v}\over{\partial\,y}}  \biggr) \,$

\setbox21=\hbox{$\displaystyle 
{{\partial}\over{\partial\,t}} \biggl( \rho\,\Bigl(
e+{{1}\over{2}}  u^{2}\Bigr) \biggr)\,\,+\,\,
{{\partial}\over{\partial\,x}}\biggl(\rho\,u\,\Bigl(
e+{{1}\over{2}}  u^{2}\Bigr)\,+\,p\,u\biggr)\,\,+\,\,$}
\setbox22=\hbox{$\displaystyle   \, 
\,+\,{{1}\over{y}}\,{{\partial}\over{\partial\,y}}\biggl(\rho\,v\,y\, \Bigl(
e+{{1}\over{2}}u^{2} \Bigr)\,+\,p\, v \,y \biggr)\,\,=\,\,
 {{\mu}\over{y}}\, {{\partial}\over{\partial\,y}}\Bigl(\,u\,
{{\partial\,u}\over{\partial\,y}}\,  \Bigr) + \,\,  {{k}\over{y}}\, 
{{\partial}\over{\partial\,y}} \, \Bigl( y \,  {{\partial\,T}\over{\partial\,y}} \Bigr)
\,$}
\setbox30= \vbox {\halign{#\cr \box21 \cr \box22 \cr}}
\setbox31= \hbox{ $\vcenter {\box30} $}
\setbox44=\hbox{\noindent  (6.5) $\,\,   \left\{ \box31 \right. $}  
\smallskip \noindent $ \box44 $

\smallskip \noindent  (6.6) $\qquad \displaystyle
0 \,\, \leq \,\, y \,\, \leq \,\, h \,.$ 

\bigskip \noindent $\bullet \quad$ 
The derivation of the coupled model can be conducted as in the previous sections. We first
introduce the mean value of density, momentum and internal energy as in (3.2), (3.3) and
(3.4)~:

\smallskip \noindent  (6.7) $\qquad \displaystyle
\widetilde{\rho} (t,x)\,\,=\,\,{{2}\over{h^2}}
\int_{\displaystyle 0}^{\displaystyle {h}} \,\rho(t,x,y)\,y \, {\rm d}y  $

\smallskip \noindent  (6.8) $\qquad \displaystyle
\widetilde{\rho} (t,x)\,\widetilde{u} (t,x)  \, \,=\,\,{{2}\over{h^2}}
\int_{\displaystyle 0}^{\displaystyle {h}} \,(\rho u)(t,x,y)\,y \,{\rm d}y  $

\smallskip \noindent  (6.9) $\qquad \displaystyle 
\widetilde{\rho} (t,x) \, \widetilde{e} (t,x) \, +\, {{1}\over{2}} \widetilde{\rho}
(t,x)\,\widetilde{u}^{2} (t,x)  \, \,=\,\,{{2}\over{h^2}}
\int_{\displaystyle 0}^{\displaystyle {h}} \,\bigl( \, \rho e\,+\,{{1}\over{2}}
\rho u^{2} \, \bigr) (t,x,y)\,y \, {\rm d}y  \,.$

\smallskip \noindent
We multiply equation (6.2) by $\,y , \,$ integrate  between $0$ and $h$ and divide by $\,
h^2~:  2 \,.$ We get 

\smallskip \noindent  $ \qquad \qquad \displaystyle
{{\partial\,\widetilde{\rho}}\over{\partial\,t}}\,\,+\,\,{{\partial}\over{\partial\,x}}
\bigl( \widetilde{\rho}\,\widetilde{u} \bigr)\,\,+\,\,{{2}\over{h^2}} \, \Bigl[ \rho 
\, v \, y \Bigr]_{ \displaystyle y=0}^{\displaystyle y=h} \,\,= \,\, 0 \,.   $

\smallskip \noindent
The term inside the brakets in the left hand side of the previous relation is null due to
the no slip boundary condition~:  

\smallskip \noindent  (6.10) $\qquad \displaystyle
u(t,\, x, \, h) \,\,= \,\, v (t,\, x, \, h) \,\,= \,\,0 , $

\smallskip \noindent 
and the conservation of mass becomes 

\smallskip \noindent  (6.11) $\qquad \displaystyle
{{\partial\,\widetilde{\rho}}\over{\partial\,t}}\,\,+\,\,{{\partial}\over{\partial\,x}}
\bigl( \widetilde{\rho}\,\widetilde{u} \bigr) \,\, = \,\, 0 \,.  $

\bigskip \noindent $\bullet \quad$ 
We make the same operation for the momentum equation (6.3)~:  

\setbox21=\hbox{$\displaystyle 
{{\partial}\over{\partial\,t}}( \widetilde{\rho} \, \widetilde{u})\,\,+\,\,
{{\partial}\over{\partial\,x}} \, \biggl( \, {{2}\over{h^2}}
\int_{\displaystyle 0}^{\displaystyle {h}} \,(\rho u^2 \,+ \,p)(t,x,y)\,y \,{\rm d}y \,
\biggr) \,+  $}
\setbox22=\hbox{$\displaystyle   \qquad \qquad  \qquad \qquad  \qquad 
+ \,\, {{2}\over{h^2}} \, \Bigl[ \rho \, y \,u \, v   \Bigr]_{ \displaystyle
y=0}^{\displaystyle y=h} \,\,= \,\, {{2 \mu}\over{h^2}} \, \Bigl[ y \, {{\partial
u}\over{\partial y}} \Bigr]_{ \displaystyle y=0}^{\displaystyle y=h} \quad .\,$}
\setbox30= \vbox {\halign{#\cr \box21 \cr \box22 \cr}}
\setbox31= \hbox{ $\vcenter {\box30} $}
\setbox44=\hbox{\noindent  (6.12) $\,\,   \left\{ \box31 \right. $}  
\smallskip \noindent $ \box44 $

\smallskip \noindent
Under the same hypotheses concerning the boundary layer presented in section 3, we have~:

\smallskip \noindent  (6.13) $\qquad \displaystyle
(\gamma -1)\,\widetilde{\rho}\,\widetilde{e} \,\,\,\simeq\,\,\, {{2}\over{h^2}}
\int_{\displaystyle 0}^{\displaystyle {h}}  \,\,p(t,x,y) \,y \, {\rm d}y  $

\smallskip \noindent  (6.14) $\qquad \displaystyle
\widetilde{\rho}\,\widetilde{u}^{2} \,\,\,\simeq\,\,\,{{2}\over{h^2}}
\int_{\displaystyle 0}^{\displaystyle {h}}  \bigl( \rho \, u^{2} \bigr)\,(t,x,y) \,y \,
{\rm d}y \,.  $

\smallskip \noindent
We insert these evaluations inside relation (6.12) and obtain 

\smallskip \noindent  (6.15) $\qquad \displaystyle
{{\partial\,\widetilde{\rho} \,  \widetilde{u} } \over{\partial\,t}}
\,\,+\,\,{{\partial}\over{\partial\,x}} \Bigl( \widetilde{\rho} \,  \widetilde{u}^2 \,+\, 
\widetilde{p}  \Bigr) \,\,= \,\, {{2 \mu}\over{h}} \, \Bigl( {{\partial u}\over{\partial
y}}\Bigr)  (t,\, x, \, h) \,$

\smallskip \noindent
with

\smallskip \noindent  (6.16) $\qquad \displaystyle
\widetilde{p}(t,x)\,\,\,\equiv \,\,\,(\gamma -1)\,\widetilde{\rho}\,\widetilde{e} \,. $

\smallskip \noindent
The treatment of the energy equation (6.5) is obtained by the same way, due to the boundary
condition (6.10), approximations (6.13), (6.14) and 

\smallskip \noindent  (6.17) $\qquad \displaystyle
\widetilde{\rho} \, \widetilde{u}\, \bigl( \widetilde{e}  \, +\, {{1}\over{2}}
\widetilde{u}^{2} \bigr) \,+\,\widetilde{p}\,\widetilde{u}
\,\,\,\simeq\,\,\,{{2}\over{h^2}} \int_{\displaystyle 0}^{\displaystyle {h}} 
\Bigl(  \rho  u \, \bigl( \, e+{{u^{2}}\over{2}} \, \bigr) +pu\, \Bigr) (t,x,y) \,y \, {\rm
d}y \,$

\smallskip \noindent  (6.18) $\qquad \displaystyle
{{\partial}\over{\partial\,t}}\bigl(\, 
\widetilde{\rho} \, \widetilde{e}  \, +\, {{1}\over{2}} \widetilde{\rho} \, 
\widetilde{u}^{2} \,\bigr)\,+\,{{\partial}\over{\partial\,x}} \biggl(
\widetilde{\rho} \, \widetilde{u}\, \bigl( \widetilde{e}  \, +\,
{{1}\over{2}} \widetilde{u}^{2} \bigr) \,+\,\widetilde{p}\,\widetilde{u}
\biggr)  \,\, = \,\,     {{2 k}\over{h}} \,{{\partial T}\over{\partial y}}\,(t,x,h)  \,\,. $

\bigskip \noindent $\bullet \quad$ 
We change the notations, replace the "tilde" unknown functions by letters without tilda,
denote by $\, \xi(t,\,x,\,\eta)\,$ the velocity in the boundary layer ($\, 0 \leq \eta <
+\infty$) and by  $\, \theta(t,\,x,\,\eta)\,$ the temperature in the same conditions. Due
to the relation 

\smallskip \noindent  (6.19) $\qquad \displaystyle
y \,\,= \,\, h - \eta \,\,, \quad \eta \approx 0 \,$ 

\smallskip \noindent
we have to change the sign in the right hand side of equations (6.15) and (6.18). We get
finally, as in (5.6)-(5.9)~:  

\smallskip \noindent  (6.20) $\qquad \displaystyle
p(t,x)\,\,\,\equiv \,\,\,(\gamma -1)\,\rho\,e$
\smallskip \noindent  (6.21) $\qquad \displaystyle
{{\partial\,\rho}\over{\partial t}}\,\,+\,\,
{{\partial}\over{\partial x}} \bigl(\, \rho \, u  \, \bigr)\,\,=\,\,0 \,$
\smallskip \noindent  (6.22) $\qquad \displaystyle
{{\partial}\over{\partial\,t}}\bigl(\, \rho \, u  \, \bigr)\,\,+\,\,
{{\partial}\over{\partial\,x}}  \bigl(\, \rho\,u^{2}\,+\, p \, \bigr) \,\,= \,\, 
-  {{2 \mu}\over{h}} \,\,{{\partial \xi}\over{\partial \eta}}\,\bigl(t,\,x,\, 0\bigr)
\, $
\smallskip \noindent  (6.23) $\qquad \displaystyle
{{\partial}\over{\partial\,t}} \bigl(\, \rho \, e  \, 
 \, +\,\,  {{1}\over{2}}\, \rho \, u^2 \,  \bigr)\,\,+\,\,
{{\partial}\over{\partial\,x}}  \bigl( \, \rho \, u \, e  \, 
 \, +\,\, {{1}\over{2}}\, \rho \, u^3 \,   \, +\,  \, p \, u \, \bigr)
\,\,= \,\,  -  {{2 k}\over{h}} \,\,{{\partial \theta}\over{\partial
\eta}}\,\bigl(t,\,x,\,0 \bigr) \,. $

\smallskip \noindent  
We remark that there is just a factor of $2$ that makes different the set of equatons
(6.22) (6.23) from the set of relations (5.8) (5.9). 

\bigskip \noindent $\bullet \quad$ 
Inside the boundary layer, section 4 can be applied in a straightforward manner. We denote
by $\, u' \,$ the (infinitesimal) velocity, by $\, p' \,$ the difference between pressure
field $ \,p \, $ and ambiant pressure $ \, p_0 \,$ and by $\, T' \,$ the difference $ \,
T-\theta_0 .\,$ We have, due to the relations (4.4) and (4.5)~:  

\smallskip \noindent  (6.24) $\qquad \displaystyle
\rho_0 \, {{\partial u'}\over{\partial t}} \,-\, {{\mu }\over{y}}\,
{{\partial}\over{\partial\,y}}\Bigl( y \,  {{\partial u'}
\over{\partial y}} \Bigr)  \,\,=\,\,-\, {{\partial p'}\over{\partial x}} $

\smallskip \noindent  (6.25) $\qquad \displaystyle
\rho_0 \, C_p \, {{\partial T'}\over{\partial t}} \,-\, {{k}\over{y}}\,
{{\partial}\over{\partial\,y}}\Bigl( y \,  {{\partial T'}
\over{\partial y}} \Bigr) \,\,=\,\, {{\partial p'}\over{\partial t}}  \,.$

\smallskip \noindent  
In the boundary layer, we have $ \, y \approx h \,$ and the curvature effects due to 
geometry are associated to the radius $h$ and this distance is very big compared with the
thickness $\, \delta \,$ of the boundary layer~:  

\smallskip \noindent  (6.26) $\qquad \displaystyle
\delta \,\, << \,\, h \, . \,$

\smallskip \noindent  
To fix the ideas, velocity field $\, u' \,$ can be expanded with an ansatz of the type 

\smallskip \noindent  (6.27) $\qquad \displaystyle
u' \,\, = \,\, U \, f\Bigl( {{y}\over{\delta}} \Bigr) \,$ 

\smallskip \noindent
where $\, f(\scriptstyle {\bullet}) \,$ is a regular function satisfying the conditions 

\smallskip \noindent  (6.28) $\qquad \displaystyle
f'(0) \,\, \approx \,\, f''(0) \,\, \approx \,\, f(1) \,\, \approx \,\, O(1) \,. $ 

\smallskip \noindent 
Then  $ \qquad \displaystyle
{{1}\over{h}} \, {{\partial u'}\over{\partial y}} \,\,  \approx \,\,  {{U}\over{h}}
\, {{1}\over{\delta}} \, f'(0) \,, \quad {{\partial ^2 u'}\over{\partial y^2}} \,\, 
\approx \,\,  {{U}\over{\delta}} \, {{1}\over{\delta}} \, f''(0) \qquad $ and due to the 
relation (6.26),  the term $\,\,  \displaystyle {{1}\over{h}} \, {{\partial
u'}\over{\partial y}} \,\,$ can be neglected in comparison with the second term
$\,\,\displaystyle  {{\partial ^2 u'}\over{\partial y^2}} . \,$

\bigskip \noindent $\bullet \quad$ 
Finally the equations in the boundary layer can be written as

\smallskip \noindent  (6.29) $\qquad \displaystyle
\rho_0 \,  {{\partial \xi}\over{\partial t}} \,\,-\,\, \mu
{{\partial^2 \xi}\over{\partial \eta^2}} \, \,\,=\,\, -\,{{\partial p} \over
{\partial x}} $

\smallskip \noindent  (6.30) $\qquad \displaystyle
\rho_0 \,C_p \,  {{\partial \theta}\over{\partial t}} \,\,-\, \, k 
{{\partial^2 \theta}\over{\partial \eta^2}} \, \,\,=\,\, {{\partial p}
\over {\partial t}} \,$

\smallskip \noindent
as in the two dimensional case. The coupled problem in the axisymmetric case is composed
by the set of equations (6.20)-(6.23) and (6.29)-(6.30).

\bigskip
\bigskip
\noindent {\smcaps 7) $\quad$ Numerical approximation of the coupled problem.} 
\smallskip \noindent $\bullet \quad$ 
The coupled system defined in section 5 is composed by five partial differential
equations (5.7)-(5.11), the state low of perfect gas (5.6), the boundary
conditions (5.3)(5.4) at the top-bottom of the pipe and by the inflow-outflow
boundary conditions (5.15)-(5.17). We discretize this system of equations in the
following manner.

\bigskip \noindent $\bullet \quad$ 
First we introduce some integer $\, J \,$ and an associated space step $\, \Delta x
\,$ : 

\smallskip \noindent  (7.1) $\qquad \displaystyle
\Delta x \,\,=\,\, {L \over J} $

\smallskip \noindent 
and some time step $\, \Delta t \,$ is chosen below. We define the discrete
variables   $\, \rho_j ,\, u_j  ,\, e_j \,$ for density, velocity and energy at
discrete point $\, x_j \,=\, j \, \Delta x \, (j = 0, 1, 2, \dots , J) \,$ and at
time $\, t^m \,=\, m \, \Delta t \,$  :

\setbox21=\hbox{$\displaystyle 
\rho_j^m \,\, \approx \,\, \rho \, (m\, \Delta t  \,,\,  j \, \Delta x ) $}
\setbox22=\hbox{$\displaystyle  
u_j^m \,\, \approx \,\, u \, (m\, \Delta t  \,,\,  j \, \Delta x ) $}
\setbox23=\hbox{$\displaystyle  
e_j^m \,\, \approx \,\, e \, (m\, \Delta t  \,,\,  j \, \Delta x ) \,$}
\setbox30= \vbox {\halign{#\cr \box21 \cr \box22 \cr  \box23 \cr  }}
\setbox31= \hbox{ $\vcenter {\box30} $}
\setbox44=\hbox{\noindent  (7.2) $\,\,   \left\{ \box31 \right. $}  
\smallskip \noindent $ \box44 $

\smallskip \noindent
and we suppose that state equation is satisfied at time step $\, m \, \Delta t
\,$ and at vertex  $\, x_j \,=\, j \, \Delta x \,$ : 

\smallskip \noindent  (7.3) $\qquad \displaystyle
p_j^m \,\,=\,\, (\gamma - 1) \, \rho_j^m \, e_j^m \,, \quad 0 \, \le \, j \, \le
\, J \,, \quad 0 \, \le \, m \, \le \,  n \,.  $

\bigskip \noindent $\bullet \quad$  
We introduce the conservative variables $\, W \,$ for mean flow : 

\smallskip \noindent $ \displaystyle
W \,\,=\,\, \Bigl( \,  \rho \,\,,\,\, \rho \,u \,\,,\,\, \rho  \, \bigl(
e+{{\displaystyle u^2}\over {\displaystyle 2}} \bigr) \, \Bigr)^{\displaystyle \rm t}
\,, $

\smallskip \noindent  
the physical flux function $\, f(W) \,$ : 

\smallskip \noindent  (7.4) $\qquad \displaystyle
 f(W) \,\,=\,\, \Bigl( \,  \rho \, u \,\,,\,\, \rho \,u^2 \,+\, p \,\,,\,\, \rho \,
u \, \bigl( e+{{\displaystyle u^2}\over {\displaystyle 2}} \bigr) \,+\, p \, u
\, \Bigr)^{\displaystyle \rm t} $

\smallskip \noindent  
and the source term due to the boundary layer : 

\smallskip \noindent  (7.5) $\qquad \displaystyle
G(W) \,\,=\,\,  - \Bigl( \, 0 \,\,,\,\, \beta \mu \, {{\partial \xi}\over
{\partial \eta}} \bigl(t,\,x,\, 0\bigr) 
\,\,,\,\, \beta k \, {{\partial \theta}\over{\partial \eta}} \bigl(t,\,x,\, 0\bigr)  \, 
\Bigr)^{\displaystyle \rm t} \,. $

\smallskip \noindent 
where the variable  $ \, \beta \, $ is defined by the condition

\setbox21=\hbox{$\displaystyle 
1 \qquad $ in the plane case of relations (5.7)-(5.9)}
\setbox22=\hbox{$\displaystyle 
2 \qquad $ in the axisymmetric  case of relations (6.21)-(6.23).  }
\setbox30= \vbox {\halign{#\cr \box21 \cr \box22 \cr   }}
\setbox31= \hbox{ $\vcenter {\box30} $}
\setbox44=\hbox{\noindent  (7.6) $\quad \beta \,\, = \,\,    \left\{ \box31 \right. $}  
\smallskip \noindent $ \box44 $

\smallskip \noindent 
Then the equations (5.7)-(5.9) and (6.21)-(6.23) can be written in a more compact form :  
 
\smallskip \noindent  (7.7) $\qquad \displaystyle
{{\partial W}\over{\partial t}} \,+\, {{\partial}\over{\partial x}}\, f(W)
\,\,=\,\, G(W)  \,. $

\bigskip \noindent  $\bullet \quad$ 
Variables (7.2) (for $\, j = 1,\, 2,\, \dots, \, J-1) $ are advanced between times
$\, t^n \,=\, n\, \Delta t \,$ and $\,   t^{n+1} \,=\, (n+1) \Delta t \,$  according
to the Lax-Wendroff [LW60] numerical scheme. This scheme is founded on a second
order Taylor expansion in time of the conserved variables : 

\smallskip \noindent  (7.8) $\qquad \displaystyle
W_j^{n+1} \,\,=\,\, W_j^n \,+\, \Delta t \, \biggl( {{\partial W}\over{\partial t}}
\biggr)_j^n \,+\, {1\over2} \, \Delta t^2 \, \biggl( {{\partial^2 W}\over{\partial
t^2}} \biggr)_j^n  $

\smallskip \noindent 
that is exact up to a third order term relatively to variable $\, \Delta t\,$ which
is omitted in the numerical scheme (7.8). The first derivative in time $ \,
\displaystyle \biggl( {{\partial W} \over{\partial t}} \biggr)_j^n \,$ is directly
evaluated thanks to equation (7.7)~: 

\smallskip \noindent   $  \displaystyle
\biggl( {{\partial W}\over{\partial t}} \biggr)_j^n \,\,=\,\, \Bigl( G(W)
\Bigr)_j^n \,-\, \biggl( {{\partial f(W)}\over{\partial x}} \biggr)_j^n \,$

\smallskip \noindent 
and more precisely 

\setbox21=\hbox{$\displaystyle 
 \Bigl( 0 \,,\, \beta  \mu \, {{\partial \xi}\over {\partial \eta}}
 \bigl(t^n,\,j\, \Delta x,\, 0\bigr)   \,,\, \beta k \, {{\partial \theta}\over{\partial
\eta}}  \bigl(t^n,\,j\, \Delta x,\, 0\bigr)  \Bigr)^{\displaystyle \rm t} \, \,\,  $}
\setbox23=\hbox{$\displaystyle  
\,\,-\,\, {{1}\over{2 \, \Delta x}} \, \biggl( \, \biggl[ \Bigl( \rho u \,,\, \rho
u^2 \,+\,p \,,\, \rho u \, ( e+{{\displaystyle u^2}\over {\displaystyle 2}}) \,+\,
pu \Bigr)^{\displaystyle \rm t} \, \biggr]_{j+1}^n \, $}
\setbox24=\hbox{$\displaystyle   \qquad  \quad  
\,\,-\,\, \biggl[ \Bigl( \rho u \,,\, \rho u^2 \,+\,p \,,\, \rho u \,
(e+{{\displaystyle u^2}\over {\displaystyle 2}}) \,+\, pu \Bigr)^{\displaystyle \rm t}
\, \biggr]_{j-1}^n \, \biggr) \,.$}
\setbox30= \vbox {\halign{#\cr \box21 \cr   \box23 \cr \box24 \cr   }}
\setbox31= \hbox{ $\vcenter {\box30} $}
\setbox44=\hbox{\noindent  (7.9) $\displaystyle  \quad \biggl( {{\partial
W}\over{\partial t}} \biggr)_j^n \,\,=\,\, \left\{ \box31 \right. $}  
\smallskip \noindent $ \box44 $

\smallskip \noindent 
The second derivative in time $\displaystyle \, \biggl( {{\partial^2
W}\over{\partial t^2}} \biggr)_j^n  \,$ is obtained by a derivation of equation
(7.7) that takes into account the Schwarz property for partial derivatives :

\smallskip \noindent   $  \displaystyle
\biggl( {{\partial^2 W}\over{\partial t^2}} \biggr)_j^n \,\,=\,\, 
\Bigl( \partial_t G(W) \Bigr)_j^n   \,-\, \Biggl[ {{\partial}\over{\partial x}}\, 
\biggl({{\partial f(W)}\over{\partial t}} \biggr)  \Biggr]_j^n \,$

\smallskip \noindent 
and after discretization of the $\displaystyle \, {{\partial}\over{\partial x}} \,$
operator by finite differences, we get : 

\setbox21=\hbox{$\displaystyle 
\Bigl( \partial_t G(W) \Bigr)_j^n  \,\,-\,\, {{1}\over{\Delta x}} \, \biggl\{ \, 
\Bigl( \partial_W f(W) \Bigr) _{j+1/2}^n \, {\scriptstyle \bullet}\, 
\bigl( \partial_t W \bigr)_{j+1/2}^n $}
\setbox22=\hbox{$\displaystyle  \qquad \qquad \qquad  
\,\,-\,\, \,  \Bigl( \partial_W f(W) \Bigr) _{j-1/2}^n \, {\scriptstyle \bullet}\, 
\bigl( \partial_t W \bigr)_{j-1/2}^n \, \biggr\} \, . $}
\setbox30= \vbox {\halign{#\cr \box21 \cr   \box22 \cr   }}
\setbox31= \hbox{ $\vcenter {\box30} $}
\setbox44=\hbox{\noindent  (7.10) $\displaystyle  \quad \biggl( {{\partial^2
W}\over{\partial t^2}} \biggr)_j^n \,\,=\,\,  \left\{ \box31 \right. $}  
\smallskip \noindent $ \box44 $

\smallskip \noindent 
We use  classical expressions for the discrete operators presented in equation
(7.10)~: the time derivative of right hand side of equation (7.7) is local in space and
will be evaluated ``more above'' : 

\smallskip \noindent  (7.11) $ \,\,\, \displaystyle
\bigl( \partial_t G(W) \bigr)_{j}^n \,= \, \biggl[ \,- 
{{\partial}\over{\partial t}} \, \Bigl( 0 \,,\, \beta  \mu \, {{\partial \xi}\over
{\partial \eta}}  \bigl(t^n,\,j\, \Delta x,\, 0\bigr) \,,\, \beta k \,{{\partial
\theta}\over{\partial \eta}}  \bigl(t^n,\,j\, \Delta x,\, 0\bigr) 
\Bigr)^{\displaystyle \rm t} \,  \biggr]_{j}^n  \,,$

\smallskip \noindent 
the jacobian matrix $\, \Bigl( \partial_W f(W) \Bigr) _{j+1/2}^n \,$ at the intermediate
point $\,(j\!+\!1/2)\Delta x\,$ is evaluated thanks to a simple two-point mean value formula
: 

\setbox21=\hbox{$\displaystyle 
{1\over2} \, \Biggl\{ \Biggl[  {{ \partial \bigl( \rho u \,,\, \rho u^2 +p \,,\, 
\rho u \,  \bigl( e+{{\displaystyle u^2}\over {\displaystyle 2}} \bigr) + pu
\bigr)^{\displaystyle \rm t} }  \over  { \partial \bigl( \rho \,,\,  \rho u \,,\,  
\rho  \,  (e+{{\displaystyle u^2}\over {\displaystyle 2}} \bigr)  \bigr)
^{\displaystyle \rm t} }}  \Biggr]_j^n  \,\,+ \,\,  $}
\setbox22=\hbox{$\displaystyle  
\,\,+ \,\, \Biggl[  {{ \partial \bigl( \rho u \,,\, \rho u^2 +p \,,\, 
\rho u \,  \bigl(e+{{\displaystyle u^2}\over {\displaystyle 2}} \bigr) + pu \bigr)
^{\displaystyle \rm t} }  \over  { \partial \bigl( \rho \,,\,  \rho u \,,\,   \rho 
\,  (e+{{\displaystyle u^2}\over {\displaystyle 2}} \bigr) \bigr) ^{\displaystyle
\rm t}}} \Biggr]_{j+1}^n \Biggr\} \,, $}
\setbox30= \vbox {\halign{#\cr \box21 \cr \box22 \cr    }}
\setbox31= \hbox{ $\vcenter {\box30} $}
\setbox44=\hbox{\noindent  (7.12) $\,\, \Bigl( \partial_W f(W) \Bigr) _{j+1/2}^n
\,\,=\,\,    \left\{ \box31 \right. $}  
\smallskip \noindent $ \box44 $

\smallskip \noindent 
and the time derivative of conservative variables at intermediate point $\,(j+1/2)\Delta
x\,$ is obtained with a centered scheme : 

\setbox21=\hbox{$\displaystyle \!\!
{1\over2} \, \biggl\{ \, \Bigl( 0 \,,\, \beta  \mu \, {{\partial \xi}\over{\partial
\eta}} \bigl(t^n,\,j\, \Delta x,\, 0\bigr) \,,\, \beta k \,
{{\partial \theta}\over{\partial \eta}} \bigl(t^n,\,j\, \Delta x,\, 0\bigr) \, 
\Bigr)^{\displaystyle \rm t}  \,\,+\,\,  $}
\setbox22=\hbox{$\displaystyle   \!\!
 \Bigl( 0 ,\, \beta  \mu \, {{\partial\xi}\over{\partial \eta}} \bigl(t^n, 
(j \! + \! 1)   \Delta x,\, 0\bigr) , \, \beta k 
\, {{\partial \theta}\over{\partial \eta}}
 \bigl(t^n,\,(j \! + \! 1)\, \Delta x,\, 0\bigr)  \Bigr)^{\displaystyle \rm t} 
\biggr\}  $}
\setbox23=\hbox{$\displaystyle   \!\!
\,\,-\,\, {{1}\over{\Delta x}} \, \biggl( \, \biggl[ \Bigl( \rho u \,,\, \rho u^2
\,+\,p \,,\, \rho u \, ( e+{{\displaystyle u^2}\over {\displaystyle 2}}) \,+\, pu
\Bigr)^{\displaystyle \rm t} \, \biggr]_{j+1}^n \, $}
\setbox24=\hbox{$\displaystyle   \qquad  \quad  
\,\,-\,\, \biggl[ \Bigl( \rho u \,,\, \rho u^2 \,+\,p \,,\, \rho u \,
(e+{{\displaystyle u^2}\over {\displaystyle 2}}) \,+\, pu \Bigr)^{\displaystyle \rm t}
\, \biggr]_{j}^n \, \biggr) \,\, .$}
\setbox30= \vbox {\halign{#\cr \box21 \cr \box22 \cr  \box23 \cr \box24 \cr   }}
\setbox31= \hbox{ $\vcenter {\box30} $}
\setbox44=\hbox{\noindent  (7.13) $ \bigl(\partial_t W \bigr)_{j+1/2}^n
=  \left\{ \box31 \right. $}  
\smallskip \noindent $ \box44 $

\bigskip \noindent  $\bullet \quad$ 
The source term $\,G(W)\,$ (relation (7.5)) is simple to represent with an integral
formula, due to the simple structure of heat equations (5.10) and (5.11). We have (see {\it e.g.}
Morse and Feshbach [MF53]) : 

\smallskip \noindent  (7.14) $\quad \displaystyle
{\rm erf} \, (\chi) \,\, \equiv \,\, {{2}\over{\sqrt{\pi}}} \int_{\displaystyle 0} 
^{\displaystyle \chi} {\rm e}^{\displaystyle - \sigma^2} {\rm d}\sigma \,\,\,,$

\smallskip \noindent  (7.15) $\quad \displaystyle
\xi(t\,,\, x \,,\, \eta) \,\,=\,\, -{1\over{\rho_0}} \int_{\displaystyle 0} 
^{\displaystyle t} \, {{\partial p}\over{\partial x}}(z,\,x) \,\,  {\rm erf} \, \Biggl[
 {{\eta}\over{{\sqrt{4 {\displaystyle{{\mu}\over{\rho_0}}}  (t-z)}}}} \Biggr] {\rm d}z\,$

\smallskip \noindent  (7.16) $\quad \displaystyle
\theta(t\,,\, x \,,\, \eta) \,\,=\,\, \theta_0 \,+\, {1\over{\rho_0 \, C_p}} \,  \int_{\displaystyle 0} 
^{\displaystyle t} \, {{\partial p}\over{\partial t}}(z,\,x) \,  {\rm erf} \, \Biggl[
{{\eta}\over{{\sqrt{4 {\displaystyle{{k}\over{\rho_0 \, C_p }}}  (t-z)}}}} \Biggr] {\rm
d}z\,,$

\smallskip \noindent 
and after derivation relatively to transverse variable $\, \eta\,$ 

\smallskip \noindent  (7.17) $\quad \displaystyle
{{\partial \xi}\over{\partial \eta}} \, (t\,,\,x\,,\,0) \,\,=\,\, -{1\over{\mu}}
\int_{\displaystyle 0}  ^{\displaystyle t} \, {{\partial p}\over{\partial x}}(t-z,\,x) \,\, 
\sqrt{{{\mu}\over{\rho_0\, \pi\, z}}}\,  {\rm d}z\,$

\smallskip \noindent  (7.18) $\quad \displaystyle
{{\partial \theta}\over{\partial \eta}} \, (t\,,\,x\,,\,0)   \,\,=\,\, {1\over{k}}
\int_{\displaystyle 0}  ^{\displaystyle t} \, {{\partial p}\over{\partial t}}(t-z,\,x) \,\, 
\sqrt{{{\mu}\over{\rho_0\, C_p \, \pi\, z}}}\,  {\rm d}z\,.$

\smallskip \noindent 
In consequence, the source term $\, \Bigl( G(W) \Bigr)_j^n \,$ is numerically evaluated
according to 

\smallskip \noindent  (7.19) $\quad \displaystyle
\Bigl( G(W) \Bigr)_j^n \,\,=\,\, \pmatrix {0 \cr \displaystyle
 \beta  \sum_{m=0}^{n-1} \int_{m \Delta t}^{(m+1)\,  \Delta t} \,  {{\partial p}\over{\partial
x}}(n\,\Delta t - z ,\, x_j)\, \sqrt{{{\mu}\over{\rho_0\, \pi\, z}}}\,  {\rm d}z  \cr 
\displaystyle -\beta   \sum_{m=0}^{n-1} \int_{m \Delta t}^{(m+1)\,  \Delta t} \,  {{\partial
p}\over{\partial t}}(n\,\Delta t - z ,\, x_j)\, \sqrt{{{\mu}\over{\rho_0\, C_p \, \pi\,
z}}}  \,  {\rm d}z \cr } \,\,. $

\smallskip \noindent 
The intermediate integrals in the second line  of right hand side of (7.19) is approached
with a two-point quadrature formula relatively to the measure  $\,\displaystyle   {{{\rm
d}z}\over{\sqrt{z}}} \,$:

\smallskip \noindent  (7.20) $\quad \displaystyle
\int_{\displaystyle a}^{\displaystyle b} \, \varphi(z) \, {{{\rm d}z}\over{\sqrt{z}}} \,\,
\approx \,\, {1\over2} \, \bigl( \varphi(a) + \varphi(b) \bigr) 
\int_{\displaystyle a}^{\displaystyle b} \, {{{\rm d}z}\over{\sqrt{z}}}  \,\,=\,\, 
 \bigl( \varphi(a) + \varphi(b) \bigr)  \, {{b-a}\over{\sqrt{a}+\sqrt{b}}}  \,$

\smallskip \noindent 
and integrals in third line of   relation  (7.19) are numerically approached by a
one-point quadrature formula : 

\smallskip \noindent  (7.21) $\quad \displaystyle
\int_{\displaystyle a}^{\displaystyle b} \, \varphi(z) \, {{{\rm d}z}\over{\sqrt{z}}} \,\,
\approx \,\,  \,  \varphi \Bigl( {{a+b}\over{2}} \Bigr) 
\int_{\displaystyle a}^{\displaystyle b} \, {{{\rm d}z}\over{\sqrt{z}}}  \,\,=\,\, 2 \, 
 \varphi \Bigl( {{a+b}\over{2}} \Bigr)   \, {{b-a}\over{\sqrt{a}+\sqrt{b}}}  \,. $

\smallskip \noindent 
We deduce, due to quadrature relations (7.20)-(7.21) and elementary use of finite
differences : 

\setbox21=\hbox{$\displaystyle 
\Bigl( G_2(W) \Bigr)_j^n \,\, =\,\,  \, $}
\setbox22=\hbox{$\displaystyle \,\,\,\,  = \,
 \beta \,  \sum_{m=0}^{n-1} {{\sqrt{\Delta t}}\over
{2 \Delta x}} \,  \Bigl[ \bigl( p_{j+1}^{n\!-\!m\!-\!1} +  p_{j+1}^{n\!-\!m} \bigr) -
\bigl(   p_{j-1}^{n\!-\!m\!-\!1} +  p_{j-1}^{n\!-\!m} \bigr) \Bigr] \,
\sqrt{{\mu}\over{\rho_0 \, \pi}} \, {1\over{\sqrt{m}+\sqrt{m+1}}} \, $}
\setbox30= \vbox {\halign{#\cr \box21 \cr \box22 \cr     }}
\setbox31= \hbox{ $\vcenter {\box30} $}
\setbox44=\hbox{\noindent     $  \left\{ \box31 \right. $}  
\smallskip \noindent $ \box44 $

\smallskip \noindent  $  \displaystyle
\Bigl( G_3(W) \Bigr)_j^n \,\,= \, \, - \,2 \beta \, \sum_{m=0}^{n-1} {{1}\over
{\sqrt{\Delta t}}} \,  \Bigl[  p_{j}^{n-m} -  p_{j}^{n-m-1}  \Bigr] \,
\sqrt{{{\mu}\over{\rho_0\, C_p \,
\pi}}} \, {1\over{\sqrt{m}+\sqrt{m+1}}} \,$

\smallskip \noindent 
and finally : 

\smallskip \noindent  (7.22) $\qquad \displaystyle
\Bigl( G(W) \Bigr)_j^n \,\,=\,\, \biggl( 0 \,,\, \Bigl( G_2(W) \Bigr)_j^n \,,\, \Bigl(
G_3(W) \Bigr)_j^n \, \biggr)  ^{\displaystyle \rm t}   \,. $

\bigskip 
\centerline { \epsfysize=4cm    \epsfbox  {fig3.epsf} }
\smallskip  \smallskip
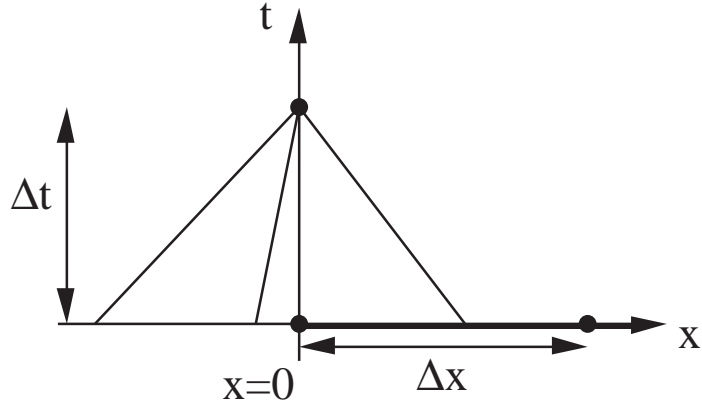
\centerline { {\bf Figure 3} \quad Characteristic directions at the entrance 
$\, x=0.\,$  }
\smallskip

\smallskip \noindent 
From previous evaluations, the time derivative of the source term $\, \bigl( \partial_t
G(W) \bigr)_{j}^n \,$  is computed with a simple first order scheme : 

\smallskip \noindent  (7.23) $\qquad \displaystyle
\bigl( \partial_t G(W) \bigr)_{j}^n \,\,= \,\,  \biggl( 0 \,,\, {{G_2(W)_j^n -
G_2(W)_j^{n-1}} \over {\Delta t}} \,,\,  {{G_3(W)_j^n - G_3(W)_j^{n-1}} \over {\Delta t}} 
\biggr) \,. \, $

\bigskip \noindent  $\bullet \quad$ 
We neglect the boundary layers when considering numerically the boundary conditions at the
input and at the output of the domain.  The boundary conditions at $\,j=0\,$ and $\,j=J\,$
are numerically implemented using the method of characteristics (see {\it e.g.} Whitham [Wh74]).
We distinguish between threen cases~: input pressure wave,  input simple velocity wave and
nonreflecting output.  In the case of an input pressure wave (at $j=0$), two
characteristics directions are going inside the computational domain (for celerities $u$
and $u+c$) and one (associated with celerity $u-c$)  is going outside (see Figure 3). We
wish to define the state $\,W_0^{n+1}\,$ at the first mesh point and at time $\,n+1\,$ ;
all the states at time level $\,n\,$ are supposed to be given and the   pressure field at
time level $\,n+1 \,$ is imposed to be equal to some numerical value  $ \, \pi^{n+1} \,$
due to the boundary condition. We denote by $\,c_0 ,\, p_0 \,$ and $\,S_0\,$ respectively
the sound celerity, the pressure and the entropy of the air at rest at usual conditions of
temperature and pressure. We first determine an external sound celerity $\,c_e \,$ and an
external velocity $\,u_e \,$ associated with a $\, C_+ \,$ input wave ; we have
classically from locally linearized theory  [Wh74]~:  

\smallskip \noindent  (7.24) $\qquad \displaystyle
u_e \,\,=\,\, {{\pi^{n+1} - p_0}\over{\rho_0 \, c_0 }}  \,$

\smallskip \noindent  (7.25) $\qquad \displaystyle
c_e \,\,=\,\, c_0 + {{\gamma-1}\over2} \, u_e \,. $

\smallskip \noindent
Secondly we interpolate data at time level $\,n\,$ and at the foot-point $P$ going backward
along the $u-c$ characteristics starting at $\, t^{n+1} \,$ from $x=0$ :  

\smallskip \noindent  (7.26) $\qquad \displaystyle
W_P \,\,=\,\, \Bigl( 1-{{u_0^n - c_0^n}\over{\Delta t}} \Bigr) \, W_0^n \,\,+ \,\, 
{{u_0^n - c_0^n}\over{\Delta t}} \, W_1^n  \, . $

\smallskip \noindent
The state $\,W_0^{n+1}\,$ is finally defined by the following three conditions~: the
characteristic variable associated to the $u-c$ wave is constant between states $\, W_P \,$
and  $\,W_0^{n+1} \, $~: 

\smallskip \noindent  (7.27) $\qquad \displaystyle
u_0^{n+1} - {{2 c_0^{n+1}} \over{\gamma -1}} \,\,= \,\, u_P - {{2 c_P} \over{\gamma -1}}
\,,  $

\smallskip \noindent
the entropy of state  $\, W_P \,$ is equal to the entropy at rest~: 

\smallskip \noindent  (7.28) $\qquad \displaystyle
S_0^{n+1} \,\,= \,\,  S_0 \,  $

\smallskip \noindent
and the characteristic variable associated to the $u+c$ wave is constant between external
state and state $\,W_0^{n+1} \, $~: 

\smallskip \noindent  (7.29) $\qquad \displaystyle
u_0^{n+1} + {{2 c_0^{n+1}} \over{\gamma -1}} \,\,= \,\, u_e + {{2 c_e} \over{\gamma -1}}
\,.  $

\smallskip \noindent
With this implementation, the single pressure datum variable $\, \pi^{n+1} \,$ is used for
two incoming waves and the outgoing wave is not reflected. We remark that, as in [DF89], 
nothing in what we have done imposes strongly the condition $\, p \bigl( W_0^{n+1} \bigr) =
\pi^{n+1} $. In some sense, this boundary condition is transparent to the outgoing waves. 

\bigskip \noindent  $\bullet \quad$ 
We use the same notations for the input simple velocity wave associated to datum $\,
U^{n+1} \,$ at the time level under study. This datum is supposed to be sufficiently small
in order to be considered as an acoustic velocity.  We first determine the celerity of an
external state $\, W_e \,$ by a relation similar to (7.24)~:

\smallskip \noindent  (7.30) $\qquad \displaystyle
c_e \,\,= \,\, c_0 + {{\gamma -1} \over{2}} \, U^{n+1} \,, $

\smallskip \noindent
we interpolate a state $\, W_P \,$ at the foot of the $u-c$ characteristic direction using
relation (7.26) and the boundary state   $\,W_0^{n+1} \, $ is computed according to
relations (7.27) along the outgoing characteristic, (7.28) along the $u$ characteristic
direction and the following relation along the $u+c$ incoming characteristic~: 

\smallskip \noindent  (7.31) $\qquad \displaystyle
u_0^{n+1} + {{2 c_0^{n+1}} \over{\gamma -1}} \,\,= \,\, U^{n+1} + {{2 c_e} \over{\gamma -1}}
\,.  $

\bigskip 
\centerline { \epsfysize=4cm    \epsfbox  {fig4.epsf} }
\smallskip  \smallskip
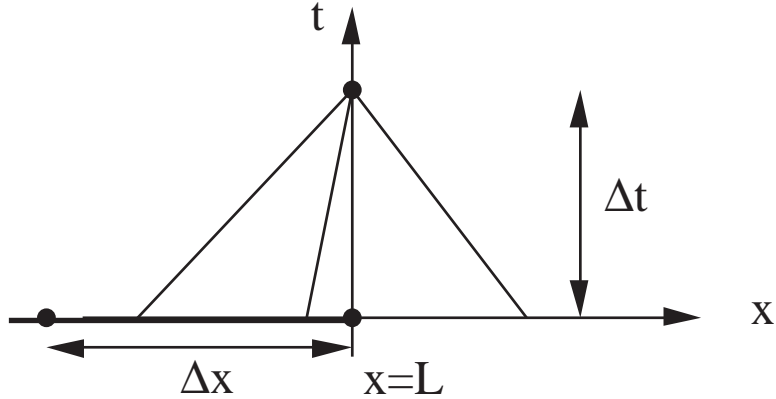
\centerline { {\bf Figure 4}	\quad 	Nonreflecting output at $\, x=L .\,$ }
\smallskip

\bigskip \noindent  $\bullet \quad$ 
For a  nonreflecting output at $x=L$ and $j=J$ (see Figure 4), the external state is the
air at rest and is  obtained by going  backward along the $u-c$ characteristic direction~: 

\smallskip \noindent  (7.32) $\qquad \displaystyle
u_J^{n+1} - {{2 c_J^{n+1}} \over{\gamma -1}} \,\,= \,\,  - {{2 c_0} \over{\gamma -1}}
\, .  $

\smallskip \noindent
Concerning the waves going outside the computational domain, we define the foot-point $Q$
associated with the $u+c$ characteristic direction with the same idea than previously~:

\smallskip \noindent  (7.33) $\qquad \displaystyle
W_Q \,\,=\,\, \Bigl( 1-{{u_0^n + c_0^n}\over{\Delta t}} \Bigr) \, W_{J-1}^n \,\,+ \,\, 
{{u_0^n + c_0^n}\over{\Delta t}} \, W_J^n  \,  $

\smallskip \noindent
and we say that the associated characteristic variable is constant between this state $W_Q$
and state   $W_J^{n+1} \,$~: 

\smallskip \noindent  (7.34) $\qquad \displaystyle
u_J^{n+1} + {{2 c_J^{n+1}} \over{\gamma -1}} \,\,= \,\, u_Q + {{2 c_Q} \over{\gamma -1}}
\,.  $

\smallskip \noindent
We suppose also than a relation similar to (7.28) determines the entropy at the limiting
vertex~: 

\smallskip \noindent  (7.35) $\qquad \displaystyle
S_J^{n+1} \,\, = \,\, S_0 \,.  $

\bigskip
\bigskip
\noindent {\smcaps 8) $\quad$  First test cases. }

\smallskip \noindent  $ \quad$   {\bf 8.1) $\quad$  Nonlinear perfect oscillating
fluid. } 

\smallskip \noindent $\bullet \quad$
In this sub-section, we neglect all the viscous effects. The continuous model is given by
equations (5.6)-(5.9) with $\, \mu = 0 \,$ and $\, k=0 \,$ and all the discrete equations
correspond to Lax-Wendroff scheme (7.8)-(7.13)  without source terms.  The first test case
consists of a simple wave going inside the domain $\, \{  x \geq 0 \}.\,$ At time $\,t \,$ 
equal to zero, the fluid is at rest (with pressure $\,p_0\,$ and temperature $\, T_0\,$
that correspond to usual thermodynamics conditions) and at $\, x=0 ,\,$ a source of
velocity $\, u(0,\,t) \,$ is supposed to be given. It defines a simple wave, submitted to
hypothesis 

\smallskip \noindent  (8.1) $\qquad \displaystyle
u \,- \, {{2 \, c} \over {\gamma - 1}} \,\,= \,\,-  {{2 \, c_0} \over {\gamma - 1}} \,$

\smallskip \noindent 
and if $\, x = X(t) \,$ is the solution of the  differential equation that defines the
characteristic line,   {\it i.e.}

\smallskip \noindent  (8.2) $\qquad \displaystyle
{{{\rm d}X} \over {{{\rm d}t}}} \,\,= \,\, u \,+\, c \, ,\, $

\smallskip \noindent
 we have (see {\it e.g.} Whitham [Wh74])

\smallskip \noindent  (8.3) $\qquad \displaystyle
u \,+ \, {{2 \, c} \over {\gamma - 1}} \,\,= \,\,{\rm Cste} \,. $

\smallskip \noindent
By elimination of sound celerity $c$ between equations (8.1) and (8.3), velocity 
$\, u({\scriptstyle \bullet}) \,$ has a constant value  along the characteristic
(8.2)-(8.3), and it is also the case for sound celerity.  We deduce that $\, u+c\,$ depends
only of its  value for $\, x = 0 \,$ and the slope of characteristic direction is
constant~: 

\smallskip \noindent  (8.4) $\qquad \displaystyle
u \,+ \, c \,\, =\,\,   c_0 \,+\, {{\gamma+1}\over{2}} \, u_0(t_0) \,. \, $ 

\smallskip \noindent 
Then characteristic lines are straight lines. Moreover, if $\,x=L\,$ is some given
abscissa, the time $\, t_L - t_0 \,$  for the wave to propagate between $\, x=0\,$ at time
$\, t= t_0\,$ and $\, x=L\,$ at time $\, t= t_L\,$ is given according to the following
relation~: 

\smallskip \noindent  (8.5) $\qquad \displaystyle
t_L \,-\, t_0 \,\,= \,\, {{L} \over {u+c}} \,\,= \,\, {{L} \over  {\displaystyle c_0 \,+\,
{{\gamma+1}\over{2}} \, u_0(t_0)}} \,. \, $

\bigskip  \noindent $\bullet \quad$
We  construct velocity field at the particular station $\, x = L\,$ with regular
time steps $\tau$. We  search velocity $\, u(L,\,n\tau) \,$ according to the  relation~: 

\smallskip \noindent  (8.6) $\qquad \displaystyle
u(L,\, n\tau) \,\, = \,\, u_0(t_{0,\,n}) \,$

\smallskip \noindent
where $\, t_{0,\,n} \, $ is the solution of the following nonlinear equation~: 

\smallskip \noindent  (8.7) $\qquad \displaystyle
n\tau \,-\, t_{0,\,n} \,\,= \,\, {{L} \over {\displaystyle  c_0 \,+\, {{\gamma+1}\over{2}} 
\, u_0(t_{0,\,n})}} \,. \, $

\smallskip \noindent 
Equation (8.7) is solved with the help of a Newton algorithm detailed in [Ms98] for a
sinusoidal input velocity $\, u_0(\theta)\,~: $ 

\smallskip \noindent  (8.8) $\qquad \displaystyle
u_0(\theta) \,\, = \,\, U_0 \, {\rm sin} \, (\omega_0 \, \theta ) \,$

\smallskip \noindent 
and the Newton algorithm is congergent without any problem as long as the characteristic
lines does not intersect,  {\it i.e.} under the condition

\smallskip \noindent  (8.9) $\qquad \displaystyle
L \,\, < \, \, L_{\rm shock} \,\, = \,\, {{2 \,  c_0^2 } \over {\displaystyle
(\gamma+1) \,\omega_0 \, U_0 }}  \, \, . \,  $

\smallskip \noindent
We introduce adimensionalized abscissa $ \, s \,$ relatively to $\, L_{\rm shock} \,$~: 

\smallskip \noindent  (8.10) $\qquad \displaystyle
s \,\, = \,\, {{x}\over{L_{\rm shock}}} \,. \, $ 

\smallskip \noindent 
The output velocity $\, u(L,\, {\scriptstyle \bullet})\,$ is a periodic function of time
with period  $\, T = 2 \pi / \omega_0 \,$ and parameter $\tau$ has been chosen such that 

\smallskip \noindent  (8.11) $\qquad \displaystyle
\tau \,\,= \,\, {{1}\over{2^N}} \, {{2 \, \pi } \over {\omega_0 }} \,\, \equiv \,\, 
{{1}\over{2^N}} \,T_{0} \, $ 

\smallskip \noindent
with a  big integer $N$ (of the order of $10$  typically) in order to proceed
a precise signal treatment. The nonlinear distorsion effect induces harmonics $\, k
\omega_0 \,$ $ \, (k=2,\,3,\, \cdots) \,$ of fondamental pulsation $\, \omega_0\,$ and
they are predicted up to $\, k=40 .\,$ Note that the time step $\tau\,$ for computing
exact solution has been chosen sufficiently small in order to avoid aliasing effects when
computing the fast Fourier transform.

\bigskip 
\centerline { \epsfysize=6cm    \epsfbox  {fig5.epsf} }
\smallskip  \smallskip
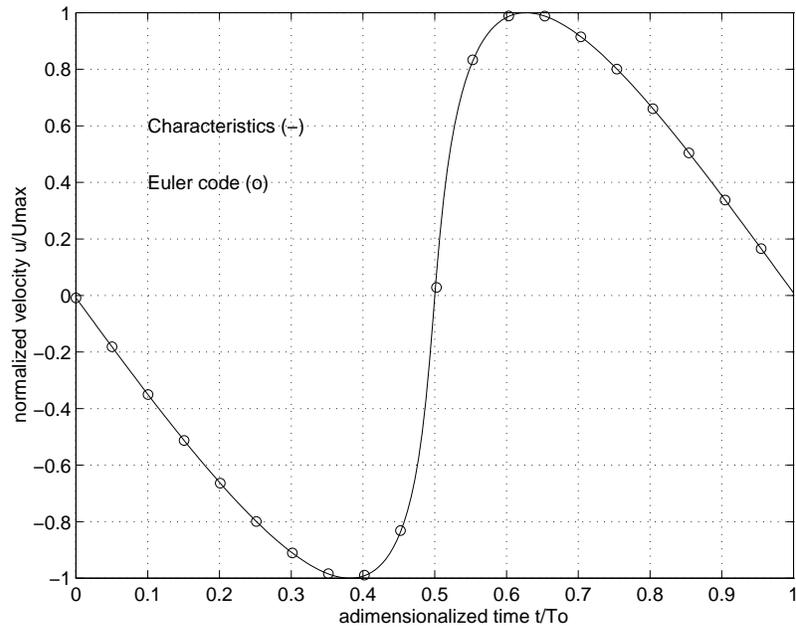
\centerline { {\bf Figure 5}	\quad 	  Signal at abscissa $\,s= 0.8$   }
\smallskip

\bigskip 
\centerline { \epsfysize=6cm    \epsfbox  {fig6.epsf} }
\smallskip  \smallskip
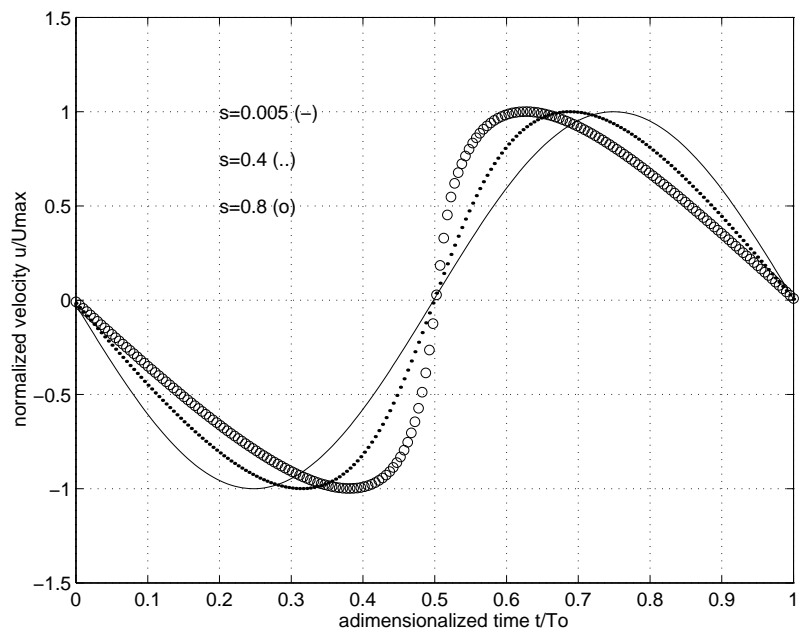
\centerline { {\bf Figure 6}	\quad 	  Signal for three abscissae    }
\smallskip

\bigskip 
\centerline { \epsfysize=6cm    \epsfbox  {fig7.epsf} }
\smallskip  \smallskip
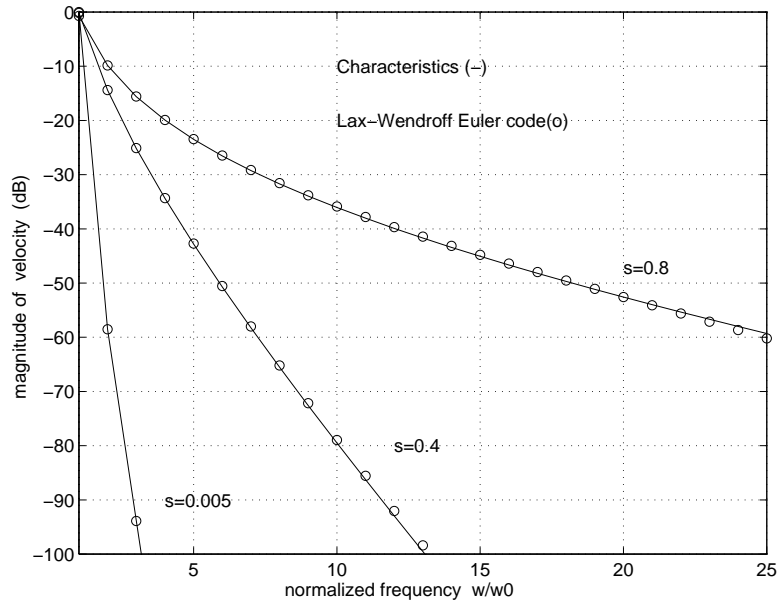
\centerline { {\bf Figure 7}	\quad Magnitude spectrum for three abscissae }
\smallskip

\bigskip 
\centerline { \epsfysize=6cm    \epsfbox  {fig8.epsf} }
\smallskip  \smallskip
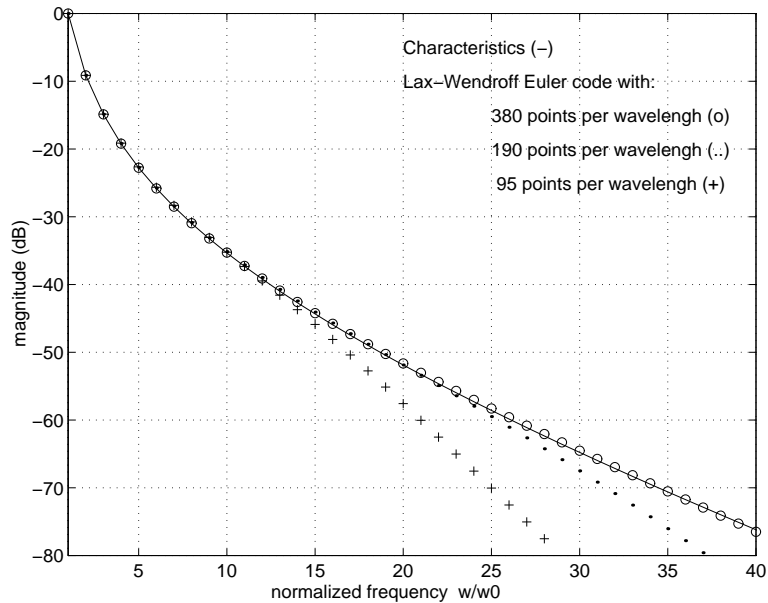
\centerline { {\bf Figure 8}	\quad 	  Influence of the space step $\, 
\Delta x $   }
\smallskip

\bigskip  \noindent $\bullet \quad$
We first test the effect of non-absorbing boundary condition on the numerical
flow computed with help of Lax-Wendroff  scheme  inside the domain. We consider two
simulations on two computational domains with the same space step $\, \Delta x .\,$ One
domain is of lenght $L$ and the other one is constructed in order to be sure that the
boundary scheme is {\bf not} active at $\, x=L .\,$ Then we  check that nonlinear
treatment (7.32)-(7.35) induces,  for waves that compose the distorted signal at $\, x=L
,\,$ relative errors in velocity that are inferior to 2 \%  for waves containing more than
10 grid points.

\bigskip  \noindent $\bullet \quad$
We  compare this simple wave with the numerical solution computed with Lax-Wendroff
scheme. We  choose a simple sinusoidal velocity profile (8.8) at the inflow (see also (7.30)
and (7.31) for the complementary boundary conditions at the inflow) and a non-reflecting
boundary condition at $\,x=L\,$ (see also relations (7.31) to (7.35)). On Figure 5, we 
plot temporal output signal for exact (characteristics) and approached (Lax Wendroff)
methods at station $\, s=0.8. \,$ We notice  that Lax-Wendroff scheme is correct for
prediction of this kind of nonlinear wave. We recover the distorsion of the wave with a
profile more and more sharp as variable $ \, s \,$ is increasing.  We  compare three
results for the Lax Wendroff scheme at $\, s = 0.005, \, 0.4 \,$ and $\, 0.8 \,$ on Figure
6 and the associated spectra for these three locations on Figure~7. We verify on Figure 7
that distorsion induces an enrichment of spectrum with a transfer of energy from low
frequency to higher frequencies. Moreover, comparison of spectra for both methods shows
that Lax-Wendroff scheme is operational for good prediction of output signal.  We compare
also spectra of output signals for different values of space steps $\, \Delta x \,$ with
constant CFL number that induces proportional values of $\, \Delta  t .\,$ For this
particular simulation ($\,u_0({\scriptstyle \bullet}) \,$ given by relation (8.8) and $\, s
= 0.8 $) we  observe (Figure 8) that numerical damping is compatible with harmonic
$\,k=15\,$ for 100 temporal points by period ({\it i.e.} $  \,  {95\over15}  \approx 6.3 $
points for this particular harmonic) and with harmonic $\,k=25 \,$ for 190 points by time
period ($ {190\over25}  = 7.6 $ points for one period).

\bigskip \noindent  $ \quad$   {\bf 8.2) $\quad$  Linear wave with visco-thermal boundary
layer effects.   } 

\smallskip \noindent $\bullet \quad$
In this sub-section, we compare our numerical model with the linear Kirchhoff theory
obtained by linearizing convective effects around a null velocity. We refer to  Bruneau et
al [BHKP89] for this classical approach in the context of first order theory with thin
boundary layer. Recall that  a wave with pulsation $ \, \omega \,$ can be a particular
solution of linear Kirchhoff theory if the phase 

\smallskip \noindent  (8.12) $\qquad \displaystyle
\Phi \,\, = \,\, \omega \, t  \,- \, K \, x \,$

\smallskip \noindent
admits a dispersion relation of the type 

\smallskip \noindent  (8.13) $\qquad \displaystyle
K \,\, = \,\, {{\omega} \over {c'(\omega) }} \,- \, i \, \alpha (\omega) \,$

\smallskip \noindent 
with

\smallskip \noindent  (8.14) $\qquad \displaystyle
c'(\omega) \,\, = \,\, c_0 \, \Biggl[\, 1 - \biggl( \sqrt{{{\mu}\over{\rho_0 \, c_0}}} \,+\,
(\gamma-1) \, \sqrt{{{k}\over{\rho_0 \, c_0 \, C_p }}}  \biggr) \, {{c_0}\over {2 \, h \,
\omega}} \,  \Biggr] \,$

\smallskip \noindent  (8.15) $\qquad \displaystyle
\alpha (\omega)\,\, = \,\,  \biggl( \sqrt{{{\mu}\over{\rho_0 \, c_0}}} \,+\,
(\gamma-1) \, \sqrt{{{k}\over{\rho_0 \, c_0 \, C_p }}}  \biggr) \,
{{c_0}\over {2 \, h \, \omega}} \, .\,$

\smallskip \noindent 
As long as the wave propagate, there is dispersion and damping of this wave. Dispersion is
due to the fact that local phase velocity $\, c'(\omega)\,$  depends on frequency (see
relation (8.14)). Damping is associated to the real part of the constant of propagation and
$\, \alpha (\omega)\,$ is a damping coefficient. 

\bigskip 
\centerline { \epsfysize=6cm    \epsfbox  {fig9.epsf} }
\smallskip  \smallskip
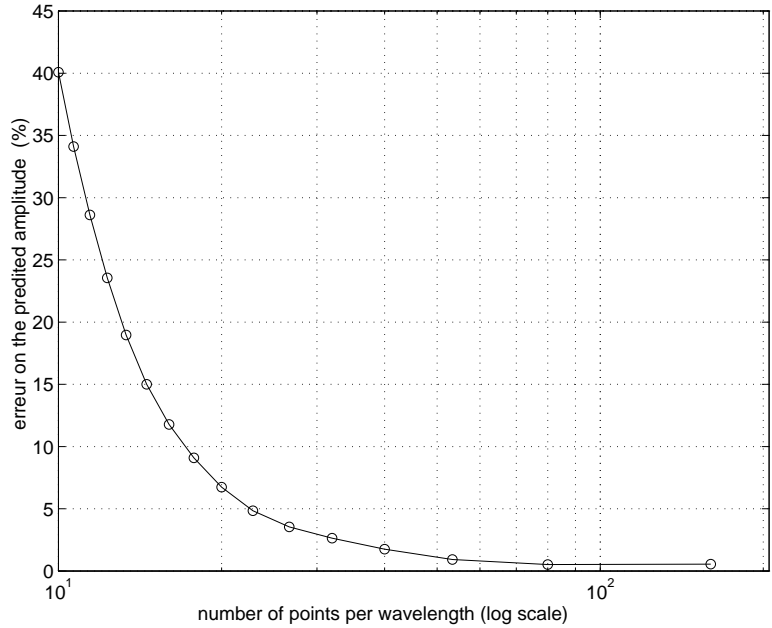
\centerline { {\bf Figure 9}	\quad 	Relative error on amplitude between  }
\centerline { the Kirchhoff theory and the Lax-Wendroff code  }
\smallskip

\bigskip  \noindent $\bullet \quad$
Sinusoidal wave of small amplitude have been simulated by the numerical model. In order to
avoid the essential of nonlinear effects, a very small amplitude  has been chosen for the
wave. We simulate the propagation of an acoustic wave.  The sinusoidal input profile is
computed over a distance $\,L \,$ physically of the order of $\, 1 \,$ meter in a pipe of
diameter of the order of $\, 1 \,$ centimeter.  We compare both amplitudes of the waves in
the first case by using  Kirchhoff model and in the second case with the pure numerical
solver. The relative errors of the predicted amplitude are plotted on Figure 9. We use as
space variable the number of grid points for one lenght wave. With more than $\, 25 \,$
points by lenght wave, we observe that relative error for pressure field  is inferior to
$\, 5 \% . \,$ The results concerning the phase obay to the same conclusion.

\bigskip 
\centerline { \epsfysize=6cm    \epsfbox  {fig10.epsf} }
\smallskip  \smallskip
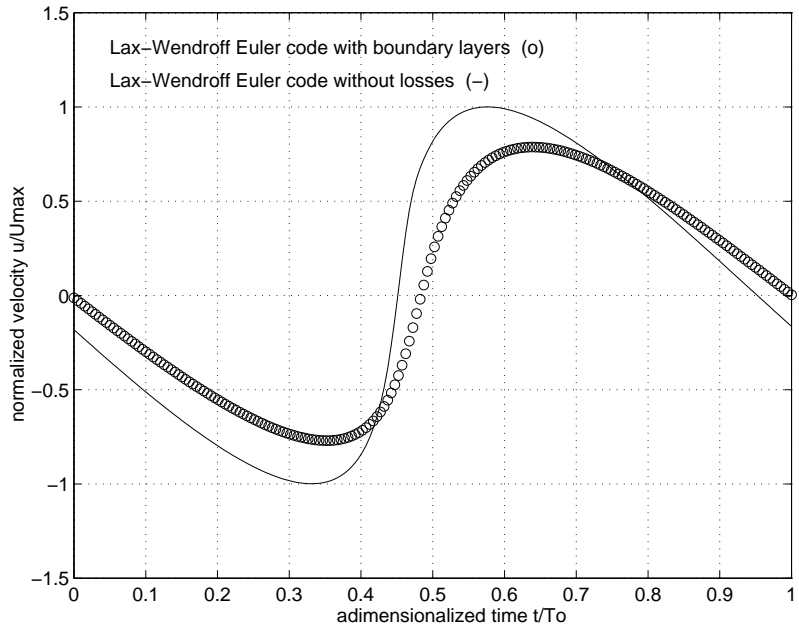
\centerline { {\bf Figure 10}	\quad 	 Comparison of flow field at 
$\, s=0.8\,$ with and  without loses   }
\smallskip

\bigskip \noindent  $ \quad$   {\bf 8.3) $\quad$  Combined nonlinear propagation and linear
boundary layer.   } 

\smallskip \noindent $\bullet \quad$
In this sub-section, we compare the shape of the wave with and without visco-thermal
boundary layer. Observe that view results are available in the literature for this kind
of elementary coupled problem. First comparison have been done with results obtained
independently by Menguy and Gilbert [MG97b]. Figure 10 presents a temporal signal of
velocity at fixed abscissa $\, s = 0.8 .\,$ There is an important damping of the wave and
we recover that wavefront is less sharp with the presence of the boundary layer. Notice
here the important remark  that viscosity associated to the boundary layer is much more
important that the one due to the thin layer Navier Stokes equations.

\bigskip \noindent  $ \quad$   {\bf 8.4) $\quad$  Trombone modelling.   } 

\smallskip \noindent $\bullet \quad $
In [HGMW96], Hirschberg,   Gilbert,  Wijnands and the first author have demonstrated
experimentally that for high level of amplitude ({\it forte, fortissimo}), there are
important nonlinear propagation effects in the trombone which can lead to shock waves. From
a modelling point of view, the slide of a trombone can be viewed in first approximation as
unidimensional pipe of lenght of the order of $\,1.5 \,$ meter and redius
$\, h = 7 $ mm. A typical incident pressure wave at the entrance of the slide is
propagating along the slide. It is a low frequency signal. At the output, we use an
absorbing boundary condition. In fact, a complete model of trombone would include the
discretization of the bell. But in the strong flairing part of the bell, the flow  is no
more quasi-unidimensional and our model is no more relevant  (see Amir, Pagneux and
Kergomard [APK97]). Because the slide is the largest lenght with a cylindrical shape of the
instrument, we conjecture that the essential of nonlinear effects occur in this part of the
instrument. We restrict our simulations to the slide alone and make the hypothesis that
nonlinear interaction between incident and reflected waves are negligeable. This hypothesis
has been verified with numerical tests [Ms98].

\bigskip 
\centerline { \epsfysize=6cm    \epsfbox  {fig11.epsf} }
\smallskip  \smallskip
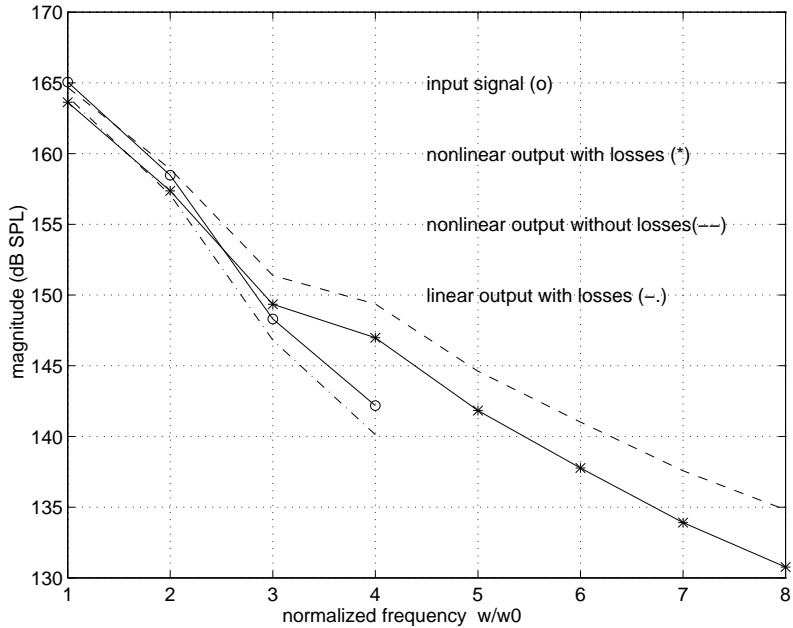
\centerline { {\bf Figure 11}	\quad 	Magnitude spectrum at the output of the slide  }
\smallskip

\bigskip \noindent $\bullet \quad $
We have synthetised a typical input signal
containing four harmonics. This signal has been propagated with a   linear propagation with
losses and with nonlinear advection with and without losses. Figure 11 represents these
three output signals.  There is of course no creation of superior harmonics with linear
propagation (done with linearized computer software developed by Quinnez [Qu95]). We
recover the four initial input modes and damping is visible (2 dB) on the sound pressure
level in Figure 11. With  nonlinear propagation without losses, new harmonics for $\, k\geq
5 \,$ are created. Moreover the amplitude of all the harmonics (except the first one) is
amplified by nonlinear propagation. This corresponds to transfer of energy from lowest
frequency to higher frequencies in order to go towards thermostatics equilibrium where we
have equi-partition of the energy between all the modes. With both nonlinear and linear
effects, previous results are damped with an amplitude varying between 1 and 5 dB for the
8 first modes. The effect of losses compensates the one of nonlinear propagation.
Nevertheless, nonlinear effects dominate the dynamics. For example the amplitude of the
4th mode is increased of 5 dB compared to the input signal. 

\bigskip \noindent $\bullet \quad $
Our numerical results confirm previous experiments done in Eindhoven~: trombone's radiated
sound is enriched by nonlinear effects which occurs in the instrument.   This effect is
known by the musicians as the "brassy" sound, typical for the trombone at loud tones (see
{\it e.g.} [GM98]). Some sonor amplitude examples are available on the net at the
following  http://www.icp.inpg.fr/$\sim$pelorson/sons.html. 


\bigskip
\bigskip
\noindent {\smcaps 9) $\quad$  Conclusion, acknowlegments.}
\smallskip \noindent $\bullet \quad $
In this study, we have proposed to use the Thin Layer Navier Stokes equations as primitive
ones to study propagation effects in thin pipes. This complete model  neglects diffusive
effects in the stream direction. Second, we have derived from this primitive set of partial
differential equations a coupled model of five equations  that takes into account both
nonlinear uni-dimensional propagation and linear diffusion in acoustic linear boundary
layer. The numerical  coupling of this two models have been done and the approach is
original~: there is no explicit need of the displacement thickness  but a set of two
velocities and two temperatures (one in the main nonviscous flow and one in the boundary
layer) allow this coupling. First numerical experiments have shown global coherence with
previous classical models in computational acoustics (characteristics, linear Kirchhoff
theory). Moreover, first application to one-dimensional modelling of trombone confirm the
importance of nonlinear wave propagation and in particular the  "brassy" sound that is
familiar to jazz musicians. The extensions of this work concern a complete treatment of
nonlinear waves with precise simulation of shock waves in the trombone,  
new numerical experiments where convolution effects in the boundary layer are computed
with direct numerical resolution of heat equation,   coupling for aerodynamic flows where
displacement effects play an important role (see {\it e.g.}  Lagree [La2k]), and mathematical
study of simplified models. The  authors thank P.Y. Lagree for  precise reading and comments 
on the first version of the manuscript.

\bigskip
\bigskip
  
\noindent {\smcaps 10) $\quad$  References.}

\smallskip \hangindent=13mm \hangafter=1 \noindent 
 [ABC92]   B. Aupoix, J. Ph. Brazier and J.
Cousteix.  Asymptotic Defect Bounda-ry-Layer Theory Applied to
Hypersonic Flows, {\it AIAA Journal} , vol.~30,  n$^{\rm o}$5, p.~1252-1259, 1992.

\smallskip \hangindent=13mm \hangafter=1 \noindent 
 [APK97]  N. Amir, V. Pagneux, J. Kergomard. Wave
propagation in acoustic horns through modal decomposition, in {\it Proceedings of the
Institut of Acoustics ISMA'97}, Edimbourgh, 1997.

\smallskip \hangindent=13mm \hangafter=1 \noindent  
[Ba67]     G.K. Batchelor.  {\it An introduction to
fluid dynamics}, Cambridge University Press, 1967.

\smallskip \hangindent=13mm \hangafter=1 \noindent   
 [BHKP89]  M. Bruneau, P. Herzog, J.
Kergomard, J.D. Polak. General formulation of the dispersion equation in bounded
visco-thermal fluid ; application to simple geometries,  {\it Wave motion}, vol.~11,
p.~441-451, 1989.

\smallskip \hangindent=13mm \hangafter=1 \noindent 
 [BL78]  B.S. Baldwin, H. Lomax. Thin layer
Approximation and Algebraic Model for Separated Turbulent Flows, AIAA Paper
n$^{\rm o}$ 78-257, AIAA 16th Aero-space Sciences Meeting, Huntsville, Alabama, 1978.

\smallskip \hangindent=13mm \hangafter=1 \noindent 
  [Bl85]  D.T. Blackstock.  Generalized Burgers equation for plane waves, 
 {\it Journal Acoust. Soc. Am.} , vol.~77, n$^{\rm o}$ 6, p.
2050-2053, 1985.

\smallskip \hangindent=13mm \hangafter=1 \noindent 
 [Br97]    Y. Brenier. Personal communication, april 1997.

\smallskip \hangindent=13mm \hangafter=1 \noindent  
 [Br98]   M. Bruneau. {\it Manuel d'acoustique
fondamentale},  Herm\`es, Paris,  1998.

\smallskip \hangindent=13mm \hangafter=1 \noindent 
 [CF48]   R. Courant, K.O. Friedrichs.   {\it Supersonic Flow and Shock Waves.} 
Interscience Publishers Inc., New York, 1948.

\smallskip \hangindent=13mm \hangafter=1 \noindent 
 [Ch64]    W. Chester.  Resonant oscillations
in closed tubes,   {\it J. Fluid Mech.}, vol.~18, p.~44-64, 1964.

\smallskip \hangindent=13mm \hangafter=1 \noindent   
[Co88]     J. Cousteix. {\it Couche limite
laminaire}, Cepadues Editions, Toulouse, 1988. 

\smallskip \hangindent=13mm \hangafter=1 \noindent    
[DF89]   F. Dubois, P. Le Floch.   Boundary Conditions for
Nonlinear Hyperbolic Systems of Conservation Laws, 
 {\it Notes on Numerical Fluid Dynamics}
(Ball\-mann-Jeltsch Editors), vol.~24, p.~96-104, Vieweg, 1989.

\smallskip \hangindent=13mm \hangafter=1 \noindent     
[GM98]   J. Gilbert, R. Msallam. Chocs cuivr\'es, 
{\it Pour la Science}, p.~27,\br
 f\'evrier~1998. 

\smallskip \hangindent=13mm \hangafter=1 \noindent      
  [Ha98]    L. Halpern.  Personal communication, april  1998.

\smallskip \hangindent=13mm \hangafter=1 \noindent       
 [He79]   G.W. Hedstrom. Nonreflecting Boundary
Conditions for Nonlinear Hyperbolic Systems, {\it J. Comput. Physics}, 
vol.~30, n$^{\rm o}$  2, p.~222-237,  1979.

\smallskip \hangindent=13mm \hangafter=1 \noindent     
   [HGMW96]    A. Hirschberg, J. Gilbert, R.
Msallam, A.P.J. Wijnands. Shock waves in trombones, {J. Acoust. Soc. Am.}, vol.~99, 
n$^{\rm o}$  3, p.~1754-1758,  1996. 

\smallskip \hangindent=13mm \hangafter=1 \noindent      
  [KCL78]    P. Kutler, S. Chakravarthy, C.K. Lombard.
AIAA Paper n$^{\rm o}$  78-213, 1978.

\smallskip \hangindent=13mm \hangafter=1 \noindent     
[Ke81]   J. Kergomard, Acoustique musicale et
champ interne des instruments \`a vent, Th\`ese  de l'Universit\'e du Maine, Le
Mans, 1981. 

\smallskip \hangindent=13mm \hangafter=1 \noindent  
   [Kr70]  H.O. Kreiss.  Initial Boundary Value
Problems for Hyperbolic Systems,  {\it Comm. Pure Appl. Math.}, 
vol.~23, p.~277-298, 1970.

\smallskip \hangindent=13mm \hangafter=1 \noindent 
  [La2k]   P.Y. Lagr\'ee.  An inverse technique to deduce 
the elasticity of a large artery, {\it European Physical Journal, Applied Physics}, 
vol.~9, p.~153-163, 2000.

\smallskip \hangindent=13mm \hangafter=1 \noindent 
  [LB80]   J.C. Le Balleur. Calcul des
\'ecoulements \`a forte interaction visqueuse au moyen de m\'ethodes de
couplage, in Computation of viscous flow interactions, AGARD CP 291, U.S.
Air Force Academy, Colorado Spring, C.O., sept.-oct. 1980. 

\smallskip \hangindent=13mm \hangafter=1 \noindent  
 [Li78]   J. Lighthill. {\it Waves in Fluids,} Cambridge
University Press, 1978. 

\smallskip \hangindent=13mm \hangafter=1 \noindent
    [LL53]   L. Landau, E. Lifschitz. {\it
Fluid Mechanics,} Nauka, Moscow, 1953. 

\smallskip \hangindent=13mm \hangafter=1 \noindent  
[LW60]   P.D. Lax, B. Wendroff. Systems of
Conservation Laws,  {\it Comm. Pure Appl. Math.}, vol.~13, p.~217-237, 1960. 

\smallskip \hangindent=13mm \hangafter=1 \noindent  
 [MDDC97]   R. Msallam, S. Dequidt, F. Dubois, 
R. Causs\'e. Mod\`ele et simulations num\'eriques de la propagation acoustique
non-lin\'eaire dans les conduits, Congr\`es of the Soci\'et\'e Fran\c{c}aise
d'Acoustique, Marseille, april 1997. 

\smallskip \hangindent=13mm \hangafter=1 \noindent  
 [MF53]   P. Morse, H. Feshbach. {\it Methods of Theoretical
Physics}, Mc Graw Hill Company, New York, 1953.

\smallskip \hangindent=13mm \hangafter=1 \noindent  
 [MG97a]   L.  Menguy, J. Gilbert, Congr\`es of the
Soci\'et\'e Fran\c{c}aise d'Acoustique, Marseille, april 1997.

\smallskip \hangindent=13mm \hangafter=1 \noindent  
 [MG97b]  L.  Menguy, J. Gilbert, personal communication, 1997.

\smallskip \hangindent=13mm \hangafter=1 \noindent  
 [Ms98]   R. Msallam. Mod\`ele et simulations
num\'eriques de l'acoustique non lin\'eai-re dans les conduits ; application \`a 
l'\'etude des effets non lin\'eaires dans le trombone. {\it Th\`ese de
l'Universit\'e Paris 6}, december 1998.

\smallskip \hangindent=13mm \hangafter=1 \noindent  
[MO97]  S. Makarov, M. Ochmann. Nonlinear and
Thermoviscous Phenomena in Acoustics, Part II,  {\it Acta Acustica} , 
vol.~83, p.~197-222, 1997.

\smallskip \hangindent=13mm \hangafter=1 \noindent 
 [MT75]  P. Merkli, H. Thoman. Transition to turbulence
in oscillating pipe flow, {\it Journal of Fluid Mechanics} , vol.~68, p.~567, 1975. 

\smallskip \hangindent=13mm \hangafter=1 \noindent  
 [Pi81]    A.D. Pierce. {\it Acoustics. An
introduction to its physical principles and applications}, Mc Graw Hill, New York, 1981.

\smallskip \hangindent=13mm \hangafter=1 \noindent    
[Qu95]   B. Quinnez. Mod\'elisation des ph\'enom\`enes
a\'ero\'elastiques bas\'ee sur une lin\'earisation des \'equations d'Euler. Th\`ese de
doctorat, {\it Ecole Centrale de Paris}, 1995. 

\smallskip \hangindent=13mm \hangafter=1 \noindent   
 [RT92]   S.G. Rubin, J.C. Tannehill. Parabolized Reduced 
Navier-Stokes Computational Techniques, {\it   Annu. Rev. Fluid Mech.}, 
 vol.~24, p.~117-144, 1992. 

\smallskip \hangindent=13mm \hangafter=1 \noindent   
 [Sc55]  H. Schlichting. {\it Boundary-Layer
Theory}, Mac Graw Hill, New York, 1955. 

\smallskip \hangindent=13mm \hangafter=1 \noindent   
 [Su91]  N. Sugimoto. Burgers equation with
fractional derivative ; hereditary effects on nonlinear acoustic waves, {\it Journal of 
Fluid Mechanics}, vol.~225, p.~631-653, 1991. 

\smallskip \hangindent=13mm \hangafter=1 \noindent    
 [Wh74]  G.B. Whitham. {\it Linear and nonlinear
waves}, John Wiley $\&$ sons, New York, 1974.

\smallskip \hangindent=13mm \hangafter=1 \noindent    
[Ze92]  R. Zeytounian. {\it Mod\'elisation
asymptotique en m\'ecanique des fluides newtoniens}, Soci\'et\'e de
Math\'ematiques Appliqu\'ees et Industrielles, Math\'ema-tiques et leurs
applications, vol.~15, Springer Verlag, 1992.

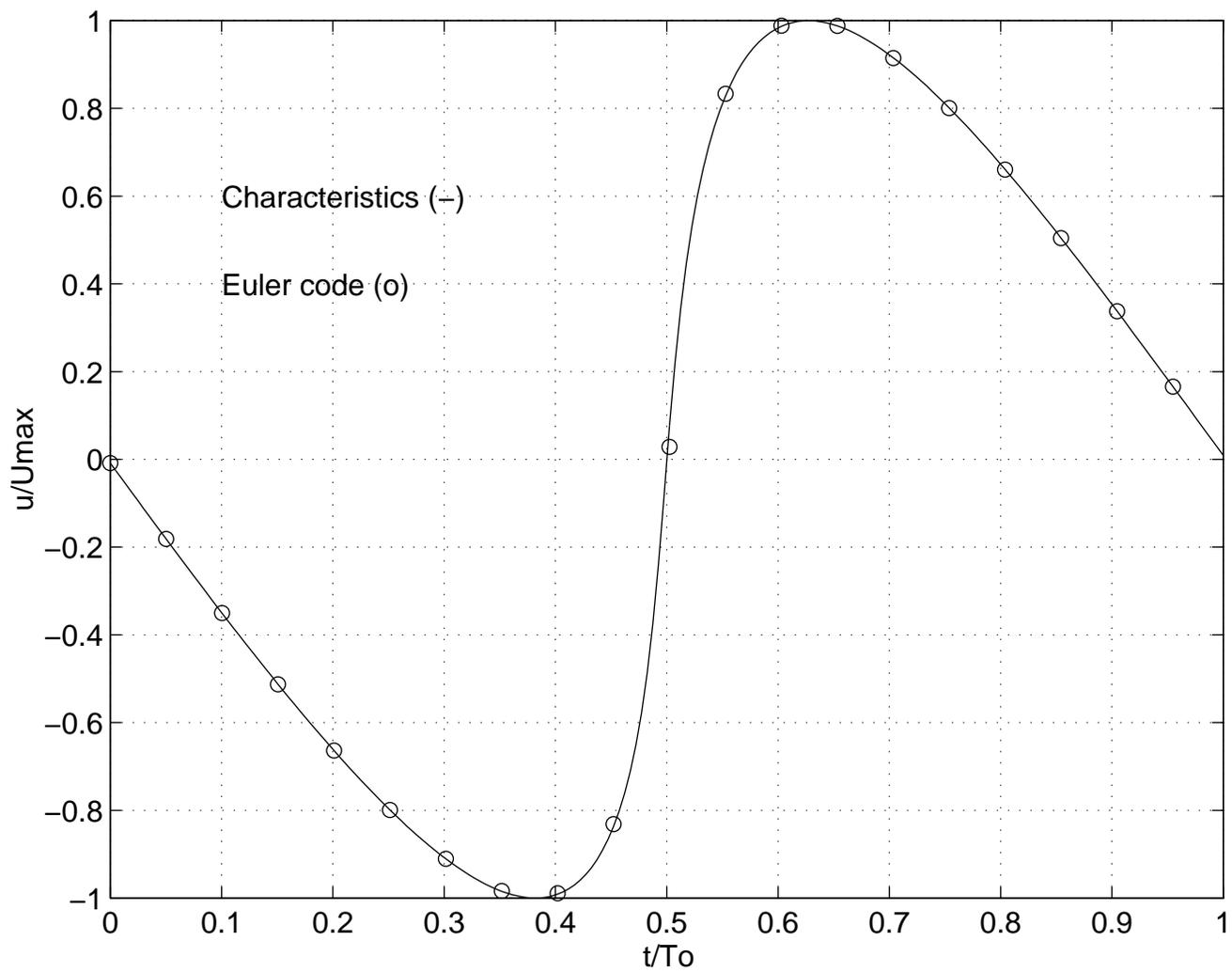
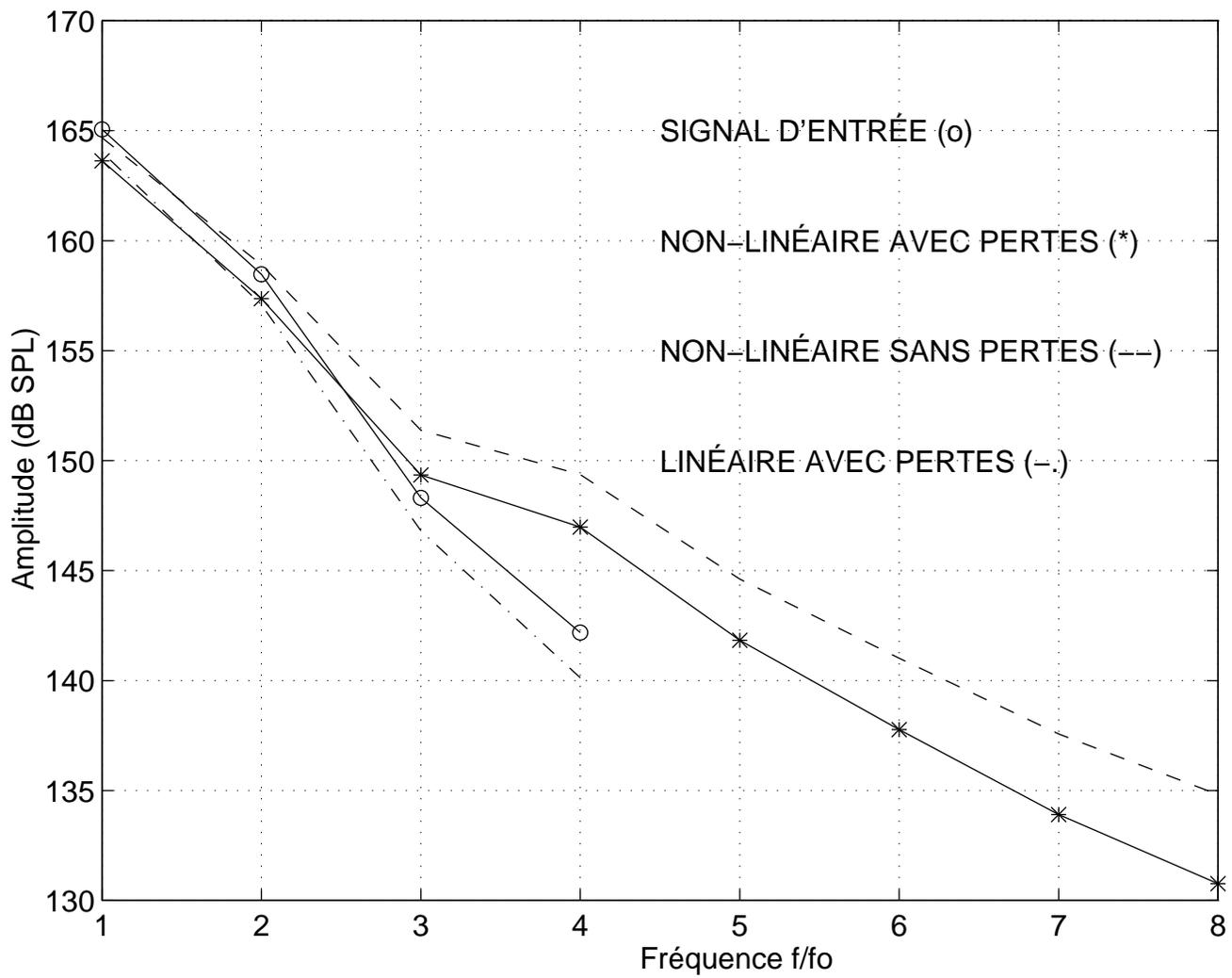
\bye